\documentclass[11pt]{article}
\usepackage[a4paper,margin=2.3cm]{geometry}
\usepackage{amsmath,amssymb,amsthm,mathrsfs}
\usepackage{xcolor}
\usepackage{parskip}
\usepackage{enumitem}
\usepackage[hidelinks]{hyperref}
\hypersetup{
  pdftitle={On L2 estimates for quadratic images of product Frostman measures},
  pdfauthor={Sung-Yi Liao, Thang Pham, Chun-Yen Shen}
}

\theoremstyle{plain}
\newtheorem{theorem}{Theorem}[section]
\newtheorem{corollary}[theorem]{Corollary}
\newtheorem{lemma}[theorem]{Lemma}
\newtheorem{proposition}[theorem]{Proposition}

\theoremstyle{definition}
\newtheorem{definition}[theorem]{Definition}

\theoremstyle{remark}
\newtheorem{remark}[theorem]{Remark}

\newcommand{\spt}{\operatorname{spt}}

\newsavebox{\strikebox}
\newcommand{\strikeword}[1]{%
  \begingroup
  \setbox\strikebox=\hbox{#1}%
  \hbox{%
    \rlap{\raisebox{.48ex}{\rule{\wd\strikebox}{0.45pt}}}%
    \box\strikebox}%
  \endgroup}

\begin{document}
\title{On $L^2$ estimates for quadratic images of product Frostman measures}
\author{Sung-Yi Liao\thanks{Department of Mathematics, Virginia Tech.
Email: \texttt{sungyi@vt.edu}}
\and Thang Pham\thanks{\raggedright Institute of Mathematics and Interdisciplinary
Sciences, Xidian University. Email: \texttt{thangpham.math@gmail.com}}
\and Chun-Yen Shen\thanks{Department of Mathematics, National Taiwan
University. Email: \texttt{cyshen@math.ntu.edu.tw}}}
\date{}

\maketitle

\begin{abstract}
Let $f\in\mathbb R[x,y,z]$ be a fixed non-degenerate quadratic
polynomial, let $0<\alpha<1$, and let $\mu$ be an
$\alpha$-Frostman probability measure supported on $[0,1]$.  We study the
pushforward measure
\[
 \nu=f_{\#}(\mu\times\mu\times\mu)
\]
on $\mathbb R$.  For every fixed nonnegative Schwartz function $\varphi$
with $\int\varphi=1$, we prove that there is $\varepsilon>0$ such that
\[
 \|\varphi_\delta*\nu\|_{L^2(\mathbb R)}^2
 \lesssim \delta^{\alpha+\varepsilon-1}
\]
for all sufficiently small $\delta>0$. As an application, if \(\mu\) is \(\alpha\)-Ahlfors--David regular on
\(A:=\spt\mu\), then
\[
 \overline{\dim}_{\mathrm B}
 \{(x,y,z)\in A^3:f(x,y,z)=t\}
 \leq 2\alpha-\varepsilon
\]
for \(\nu\)-almost every \(t\).  This is a strict improvement over the sharp general
bound \(2\alpha\), which holds uniformly in \(t\).  We also obtain a
Hausdorff-dimension bound for the exceptional set of parameters.
These conclusions give progress on a one-dimensional, three-fold
analogue of the fractal regular-value problem studied by Eswarathasan,
Iosevich, and Taylor (Advances in Mathematics, 2011).
\end{abstract}
\tableofcontents
\section{Introduction}

The study of images of fractal sets under nonlinear maps is a central theme
in geometric measure theory and harmonic analysis.

Given a map \(F\colon\mathbb R^{\ell}\to\mathbb R^m\) and a compact
set \(X\subset\mathbb R^\ell\), one seeks to relate the dimension and measure of \(F(X)\) to dimensional and
non-concentration properties of \(X\) and, when \(F\) is sufficiently
regular, to geometric features such as curvature and transversality.
A fundamental model is
the celebrated Falconer distance problem.  For a compact set \(E\subset\mathbb R^d\),
\(d\geq2\), define
\[
    \Delta(E):=\{|x-y|:x,y\in E\}.
\]
The Falconer distance conjecture asserts that
\(\mathcal L^1(\Delta(E))>0\) whenever
\(\dim_{\mathrm H}E>d/2\).

Motivated by this circle of questions, Koh, Pham, and Shen~\cite{KPS}
introduced the following definition of Falconer-type function.

\begin{definition}\label{def1.1}
Let \(\Phi: \mathbb{R}^{\ell}\)
\(\to \mathbb{R}\). For \(0<\epsilon<1\), we say that \(\Phi\) is a Falconer-type function
with exponent \(\epsilon\) if, for every compact set \(A\subset\mathbb R\)
with \(\dim_{\mathrm H}A>1-\epsilon\), the image
\(\Phi(A,\ldots,A)\) has positive Lebesgue measure.
\end{definition}

\noindent

This definition is analogous to the notion of a moderate expander over finite fields. In that setting one asks whether, for constants \(c,C>0\)
independent of \(q\), the condition
\(|A|\ge Cq^{1-\epsilon}\) forces
\(|\Phi(A,\ldots,A)|\ge cq\).

For bivariate polynomials of bounded degree over finite fields of
sufficiently large characteristic, Tao~\cite{tao} proved that moderate expansion holds unless the polynomial belongs to
one of the additive or multiplicative exceptional classes described in his
theorem.  Improved exponents are known for several
multivariable quadratic maps; see, for example,
\cite{shpas,vinh,chap,ha2,sak,benet,MP,chao}.

Koh, Pham, and Shen~\cite{KPS} proved that \textit{non-degenerate}
quadratic polynomials are Falconer-type functions with
exponent $1/3$. Here and throughout, a quadratic polynomial \(f\in\mathbb R[x,y,z]\) is called non-degenerate if $f$ depends nontrivially on each variable and is
not of the form $G(I(x)+J(y)+K(z))$ for some polynomials
\(G,I,J,K\in\mathbb R[t]\).

\begin{theorem}[\cite{KPS}, Theorem 1.2]\label{quadratic}
Let $f\in \mathbb{R}[x, y, z]$ be a non-degenerate quadratic polynomial. For a compact set $A\subset \mathbb{R}$, if \(\dim_{\mathrm H}A\)
\(>\frac{2}{3}\), then $f(A, A, A)$ has positive Lebesgue measure. 
\end{theorem}

\noindent

The result in \cite{KPS} is in fact proved for three possibly different compact
sets: if
\(\dim_{\mathrm H}A+\dim_{\mathrm H}B+\dim_{\mathrm H}C>2\), then
\(f(A,B,C)\) has positive Lebesgue measure.  This three-set dimensional condition is sharp in general.  Related results for more general
configuration maps and broader classes of functions appear in
\cite{GIT,quy}.

\noindent

In the complementary range, Pham's general theorem
\cite[Theorem~1.4(ii)]{quy} gives the following specialization to \(A=B=C\).

\begin{theorem}[\cite{quy}, Theorem 1.4(ii)]\label{quadraticsmall}
Let $f\in \mathbb{R}[x, y, z]$ be a non-degenerate quadratic polynomial. For a compact set $A\subset \mathbb{R}$, if \(\dim_{\mathrm H}A\)
\(=s\le \frac{2}{3}\), then the Hausdorff dimension of $f(A, A, A)$ is at least $\frac{3s}{2}$. 
\end{theorem}

Taken together, Theorems~\ref{quadratic} and~\ref{quadraticsmall}
provide a continuous counterpart of the prime-field expansion theorem
in~\cite{phamV}.  Indeed, the finite-field growth exponent \(3/2\)
becomes the Hausdorff-dimension lower bound \(3s/2\) in the subcritical
range \(s\leq 2/3\), while above the critical threshold \(s=2/3\), the
image has positive Lebesgue measure.

The purpose of the present paper is to go beyond these image-set results and study the distribution of the corresponding image measure. Before stating our main theorem, we recall the following definition.

\begin{definition}\label{def:AD_alpha}
Let $\alpha > 0$. We say a probability measure $\mu$ on $\mathbb{R}^d$ is $\alpha$-Frostman with the Frostman constant $C_\mu > 0$ if for every $x \in \mathbb{R}^d$ and every $r > 0$, one has
\[
\mu\bigl(B(x,r)\bigr) \le C_\mu r^{\alpha}.
\]
\end{definition}
When $C_\mu$ is an absolute constant and the support of $\mu$ is clear from the context, we simply say $\mu $ is $\alpha$-Frostman.

Let \(\varphi\) be a fixed nonnegative Schwartz function on \(\mathbb R\)
with \(\int_{\mathbb R}\varphi=1\), and put
\[
 \varphi_\delta(t):=\delta^{-1}\varphi(t/\delta).
\]
For an \(\alpha\)-Frostman probability measure \(\mu\) supported on
\([0,1]\), define
\[
 \nu:=f_\#(\mu^{\otimes3}).
\]
Equivalently,
\[
 \int h(t)\,d\nu(t)
 =\int h\bigl(f(x,y,z)\bigr)\,d\mu(x)d\mu(y)d\mu(z)
\]
for every bounded Borel function \(h\).  

Instead of asking only whether the support of \(\nu\), namely
\(f(\spt\mu,\spt\mu,\spt\mu)\),
has positive Lebesgue measure, we investigate the smoothed \(L^2\)-energy of
\(\nu\) at scale \(\delta\). This is a more quantitative problem: it measures
not merely the size of the image set, but also how evenly the mass of
\(\mu\times\mu\times\mu\) is distributed under the polynomial map \(f\).

The estimate proved in this paper gives a power saving over the natural
Frostman-scale bound. Roughly speaking, it shows that for every
\(\alpha\in(0,1)\), the smoothed image measure \(\varphi_\delta*\nu\) enjoys
nontrivial \(L^2\)-control:
\[
        \|\varphi_\delta*\nu\|_{L^2(\mathbb R)}^2
        \lesssim \delta^{\alpha+\varepsilon-1}
\]
for some \(\varepsilon>0\). The gain \(\delta^\varepsilon\) is the main point.
It reflects a quantitative expansion phenomenon: the quadratic map does not
merely enlarge the support of the measure, but also prevents excessive
concentration of the image measure at small scales.

Our main result is stated as follows.

\begin{theorem}\label{thm-main-1}
Let \(f\in\mathbb R[x,y,z]\) be a fixed non-degenerate quadratic polynomial, and let \(0<\alpha<1\).  There exists
\(\varepsilon=\varepsilon(\alpha,f)>0\) such that, if \(\mu\) is an
\(\alpha\)-Frostman probability measure supported on \([0,1]\), with
Frostman constant \(C_\mu\), then
\[
 \|\varphi_\delta*\nu\|_{L^2(\mathbb R)}^2
 \leq C\delta^{\alpha+\varepsilon-1}
 \qquad(0<\delta\leq\delta_0).
\]
Here \(C\) and \(\delta_0>0\) may depend on
\(\alpha,f,\varphi\), and \(C_\mu\), but not on \(\delta\).
\end{theorem}
We do not attempt to optimize the gain \(\varepsilon\).

The conclusion does not assert that the smoothed energy remains bounded
as \(\delta\downarrow0\).  Indeed, Proposition~\ref{prop:unboundedness} gives,
for \(0<\alpha<1/2\), a non-degenerate quadratic polynomial and an
\(\alpha\)-Ahlfors--David regular measure for which
\(\|\varphi_\delta*\nu\|_2^2\gtrsim\delta^{2\alpha-1}\to\infty\).

On the other hand, suppose that \(\mu\) is
\(\alpha\)-Ahlfors--David regular on \(A:=\spt\mu\).
Proposition~\ref{prop:level-set-uniform-bound} gives the uniform estimate
\[
 \overline{\dim}_{\mathrm B}
 \{(x,y,z)\in A^3:f(x,y,z)=t\}
 \leq 2\alpha
 \qquad\text{for every }t\in\mathbb R,
\]
and the exponent \(2\alpha\) is sharp in general. Nevertheless,
Theorem~\ref{thm:level-set-typical-exceptional} shows that, after possibly
decreasing \(\varepsilon\),
\[
 \overline{\dim}_{\mathrm B}
 \{(x,y,z)\in A^3:f(x,y,z)=t\}
 \leq 2\alpha-\varepsilon
\]
for \(\nu\)-almost every \(t\), together with a Hausdorff-dimension bound
for the exceptional set of parameters. Thus, the \(L^2\)-gain yields a
genuine dimensional drop for typical level sets below the sharp uniform
threshold.

This conclusion may be viewed as progress on a one-dimensional,
three-fold instance of the fractal regular-value problem studied by
Eswarathasan, Iosevich, and Taylor~\cite{Eiosevich}. Their
Fourier-integral-operator approach treats sufficiently
large-dimensional subsets of \(\mathbb R^d\) with \(d\geq 2\), and has been extended for \(k\)-fold level sets in \cite{gra}. In contrast, the present result applies to \(A\subset\mathbb R\) for
every \(0<\alpha<1\), a regime not covered by their results and the Fourier-integral-operator approach.

\paragraph{Main ideas and novelties.}
The proof begins by rewriting the $L^{2}$ energy in a form that exposes a geometric coincidence problem, first under the additional assumption that $\mu$ is Ahlfors-David regular.
Expanding the square, inserting the definition of $\nu=f_{\#}(\mu\times \mu\times \mu)$, and applying Fubini's theorem, the estimate is reduced to bounding a six-fold integral that measures, with respect to $\mu^{\otimes 3}\times \mu^{\otimes 3}$, the mass of approximate solutions to
\[
f(x_1,y_1,z_1)\approx f(x_2,y_2,z_2)
\]
at scale $\delta$.
Equivalently, one needs to control a $\delta$-thickened energy integral of the form
\[
\int \mathbf 1_{\{|f(x_1, y_1, z_1)-f(x_2, y_2, z_2)|\lesssim \delta\}}
\,d(\mu\times \mu\times \mu)(x_1, y_1, z_1)\,d(\mu\times \mu\times \mu)(x_2, y_2, z_2).
\]
The argument then partitions this six-fold integral into a \textit{transverse} contribution, where at least one coordinate derivative of $f$ is
$\gtrsim \delta^\kappa$, for some $\kappa>0$, and a \textit{small-gradient} contribution, where all relevant derivatives are $\lesssim \delta^\kappa$ and the interaction is confined to a $\delta^\kappa$-neighbourhood of the affine critical set $K=\{\nabla f=0\}$.

In the transverse contribution, the key input is an incidence theorem for
$\delta$-separated lines and point sets of the form
$\Phi(M_1\times M_2)$, where the two factors are $\delta$-separated,
satisfy upper cardinality and non-concentration estimates, and $\Phi$ is
bi-Lipschitz on their rectangular product with controlled
$\delta$-dependent distortion; see Theorem~\ref{our_incidence_theorem}.
The proof uses the rectangular incidence theorem of Demeter and O'Regan \cite{DemeterORegan},
together with two regularization arguments that remove the losses caused by
point multiplicity and by concentration in the incident line fibres.  This
gives the power saving needed in the six-variable coincidence estimate.

If the critical set $K=\{\nabla f=0\}$ is empty, the small-gradient region is
empty at all sufficiently small scales.  Otherwise non-degeneracy forces the
Hessian of $f$ to have rank at least two, and hence $K$ is a point or an
affine line.  The small-gradient contribution is then controlled by Frostman
tube estimates, rank-two slice bounds, and a separate estimate for the
remaining full-rank case.
A weighted dyadic decomposition and a partition into uniformly
non-concentrated classes extend the transverse estimate from regular
discretizations to arbitrary Frostman measures.

\paragraph{Structure of the paper.}
Section~2 proves the bounded-window incidence theorem for bi-Lipschitz
images of rectangular products.  Section~3 proves
Theorem~\ref{thm-main-1}, first isolating the transverse counting statement
and then treating the small-gradient contribution and the passage to general
Frostman measures. Section 4 contains level-set
consequences of the $L^2$ estimate. All auxiliary algebraic and geometric
lemmas are stated at the point where they are used.

\noindent

\section{Incidences for bi-Lipschitz images of Cartesian products}
\label{sec:bounded-product-incidences}

This section establishes the incidence estimate used in the proof of
Theorem~\ref{thm-main-1}.  Extracting a power saving requires an
incidence estimate that detects the product structure of the point set.
We formulate the bounded-window version needed here and derive it from
the rectangular incidence theorem of Demeter and O'Regan
\cite[Theorem~1.5]{DemeterORegan}, using two regularization lemmas to
handle point multiplicities and concentration in the incident line
fibres.

The proof proceeds in five steps.  We first discretize the product image
as a weighted collection of $\delta$-squares.  We then decompose the
weights into rectangular Katz--Tao point layers and split the squares
according to their tube degrees.  Next, we partition the tubes so that
each incident fibre satisfies a Katz--Tao condition in line-parameter
space.  Finally, we apply the rectangular incidence theorem of Demeter
and O'Regan to each pair of layers and balance the low- and high-degree
contributions.

The resulting formulation allows two different one-dimensional factors
and requires only upper cardinality bounds.  Both features will be
needed when the Frostman measure is decomposed into non-concentrated
classes.

We first specify the metric on the space of affine lines.  For
$(\omega,b)\in S^1\times\mathbb R$, let
\[
 \ell_{\omega,b}:=\{z\in\mathbb R^2:z\cdot\omega=b\}.
\]
The pairs $(\omega,b)$ and $(-\omega,-b)$ represent the same unoriented
line.  On the corresponding quotient space $\mathcal A$, define
\begin{equation}
 d_{\mathcal A}(\ell_{\omega,b},\ell_{\omega',b'})
 :=\min_{\varepsilon\in\{-1,1\}}
 \bigl(|\omega-\varepsilon\omega'|+|b-\varepsilon b'|\bigr).
 \label{eq:inc-line-metric}
\end{equation}
This definition is independent of the chosen representatives.  On each
compact subset of $\mathcal A$, the metric $d_{\mathcal A}$ is
bi-Lipschitz equivalent to the Euclidean metric in a finite collection of
normal-coordinate charts.  We write
\[
 \ell^{(r)}:=\{z\in\mathbb R^2:\operatorname{dist}(z,\ell)\leq r\}.
\]

\begin{theorem}\label{our_incidence_theorem}
Let $0<\alpha<1$, and put
\begin{equation}
 \sigma:=\min\{2\alpha,2-2\alpha\}.
 \label{eq:inc-sigma}
\end{equation}
Fix $R,C_0\geq1$ and $0\leq\kappa<\sigma/2$, and define
\begin{equation}
 G:=\alpha\left(\frac{\sigma}{2}-\kappa\right)>0.
 \label{eq:inc-G}
\end{equation}
For every $\varepsilon$ satisfying
\begin{equation}
 0<\varepsilon<\frac{G}{4},
 \label{eq:inc-epsilon-range}
\end{equation}
there exist $\delta_0>0$ and $C<\infty$, depending only on
$\alpha,\kappa,\varepsilon,R$, and $C_0$, such that the following holds
whenever $0<\delta\leq\delta_0$.

For $i\in\{1,2\}$, let $M_i\subset[0,1]$ be $\delta$-separated and
suppose that
\begin{equation}
 |M_i|\leq C_0\delta^{-\alpha},
 \qquad
 |M_i\cap I|\leq C_0\left(\frac{|I|}{\delta}\right)^\alpha
 \label{eq:inc-M-Frostman}
\end{equation}
for every interval $I$ with $|I|\geq\delta$.  Let
$F:[0,1]^2\to\mathbb R$, and assume that the map
\[
 \Phi(x,y):=(F(x,y),y)
\]
satisfies
\begin{equation}
 C_0^{-1}\delta^\kappa|u-v|
 \leq |\Phi(u)-\Phi(v)|
 \leq C_0\delta^{-\kappa}|u-v|
 \qquad (u,v\in M_1\times M_2).
 \label{eq:inc-Phi-bilip}
\end{equation}
Set
\[
 P:=\Phi(M_1\times M_2),
\]
and suppose that $P\subset B(0,R)$.  If $L$ is a
$\delta$-separated family in $(\mathcal A,d_{\mathcal A})$ such that
\begin{equation}
 |L|\leq C_0\delta^{-2\alpha},
 \label{eq:inc-line-cardinality}
\end{equation}
then
\begin{equation}
 I_\delta(P,L)
 :=\#\{(p,\ell)\in P\times L:p\in\ell^{(\delta)}\}
 \leq C\delta^{-3\alpha+\varepsilon}.
 \label{eq:inc-main-bound}
\end{equation}
\end{theorem}

\begin{remark}[Comparison with the general incidence theorem of Fu and Ren]
\label{rem:inc-Fu-Ren}
Comparing to Theorem \ref{our_incidence_theorem}, Fu and Ren proved a similar incidence result \cite[Theorem 1.4]{FuRen} that apply to a more general case. However, such a general incidence theorem does not by itself provide the estimate required here when $0<\alpha\leq 1/2$.
Indeed, in the balanced Katz--Tao model in which both the point family
and the line-parameter family have dimension $2\alpha$, the sharp
exponent they obtain is $3\alpha$. Consequently, their theorem gives, up to subpolynomial losses, $I_\delta\lesssim_\eta\delta^{-3\alpha-\eta}$, rather than the power-saving estimate, $I_\delta\lesssim\delta^{-3\alpha+\varepsilon}$, needed below. In other words, their sharp examples show that isotropic non-concentration alone cannot force such a saving.

The point sets considered here have the additional rectangular quasi-product structure inherited from $M_1\times M_2$. The rectangular incidence theorem of Demeter and O'Regan \cite[Theorem~1.5]{DemeterORegan}, together with the two regularization lemmas below, detects this additional structure and yields the required positive power saving.
\end{remark}

\begin{remark}\label{rem:inc-bounded-window}
Only the values of $F$ on $M_1\times M_2$ occur in the theorem.  In
particular, the bounded-window assumption follows from any
scale-independent bound for $F$ on $[0,1]^2$.  Once lines whose
$\delta$-neighbourhoods miss $P$ have been discarded, their normal
parameters satisfy $|b|\leq R+1$.  Their axial lines therefore lie in a
fixed compact subset of $\mathcal A$.  This is the point at which both the
bounded physical window and the quotient metric
\eqref{eq:inc-line-metric} enter the argument.
\end{remark}

We next record the form of the incidence theorem of Demeter and O'Regan
that will be used.  If $q$ is a $\delta$-square, let $c(q)$ denote its
centre.  A collection $\mathcal Q$ of $\delta$-squares is called
rectangular $(\delta,2s,C_Q)$-Katz--Tao if
\begin{equation}
 \#\{q\in\mathcal Q:c(q)\in B\}
 \leq C_Q\left(\frac{\sqrt{rr'}}{\delta}\right)^{2s}
 \label{eq:inc-rectangular-KT}
\end{equation}
for every rectangle $B$ of arbitrary orientation and side lengths
$r\geq r'\geq\delta$.  A $\delta$-separated line family
$\mathcal Y\subset\mathcal A$ is called
$(\delta,t,C_Y)$-Katz--Tao if
\begin{equation}
 \#\bigl(\mathcal Y\cap B_{\mathcal A}(\ell,r)\bigr)
 \leq C_Y(r/\delta)^t
 \qquad\text{ for all }r\geq\delta.
 \label{eq:inc-line-KT}
\end{equation}
Here and below, a $\delta$-tube is the $\delta$-neighbourhood of its
axial line, restricted to the fixed bounded physical window under
consideration.

\begin{theorem}[Demeter--O'Regan
{\cite[Theorem~1.5]{DemeterORegan}}, coordinate form]
\label{thm:inc-DO-rectangular}
Let $0\leq s\leq1$, and put $\sigma_s:=\min\{2s,2-2s\}$.  Let
$\mathcal Q$ be a rectangular $(\delta,2s,C_Q)$-Katz--Tao collection of
$\delta$-squares in a fixed bounded square, and let $\mathcal T$ be a
$\delta$-separated collection of $\delta$-tubes.  For each
$q\in\mathcal Q$, set
\[
 \mathcal T(q):=\{T\in\mathcal T:q\cap T\neq\varnothing\}.
\]
Assume that the axial lines of $\mathcal T(q)$ form a
$(\delta,\sigma_s,C_Y)$-Katz--Tao family for every $q\in\mathcal Q$.
Then, for every $\eta>0$,
\begin{equation}
 I(\mathcal Q,\mathcal T)
 :=\#\{(q,T)\in\mathcal Q\times\mathcal T:q\cap T\neq\varnothing\}
 \leq C_{s,\eta,C_Q,C_Y}\,
 \delta^{-\eta-2s/3}|\mathcal T|^{2/3}|\mathcal Q|^{1/3}.
 \label{eq:inc-DO-bound}
\end{equation}
The constant may also depend on the fixed physical and line-parameter
windows.  The conclusion is unchanged, up to its implicit constant, if
the square side lengths, tube widths, or separation scales are changed
by fixed multiplicative factors.
\end{theorem}

Theorem~\ref{thm:inc-DO-rectangular} is a fixed-constant normalization
of \cite[Theorem~1.5]{DemeterORegan}; throughout, we use the version of
that preprint specified in the bibliography. We record the coordinate and bounded-window reduction explicitly:
\par
Choose a dyadic scale \(\bar{\delta}\) such that
\[c\delta\leq \bar{\delta}\leq C\delta,\]
where \(c,C>0\) depend only on the fixed physical window. After applying
a fixed translation and similarity, the physical window is contained
in \([0,1]^2\). Under this transformation, \(\delta\)-squares and
\(\delta\)-tubes become objects of dimensions comparable to
\(\bar{\delta}\). Each such object can be covered by \(O(1)\) standard
\(\bar{\delta}\)-squares or \(\bar{\delta}\)-tubes, respectively.
Because the original squares and axial lines are separated at a scale
comparable to \(\delta\), these replacements have uniformly bounded
multiplicity. Thus every original incidence produces only \(O(1)\)
incidences in the standard discretized model.

The rectangular Katz--Tao hypothesis is stable under this
normalization. Indeed, a similarity multiplies the side lengths
\(r,r'\) and the discretization scale by the same fixed factor, so
\[\left(\frac{\sqrt{rr'}}{\delta}\right)^{2s}\asymp\left(\frac{\sqrt{\bar r\bar r'}}{\bar\delta}\right)^{2s}.\]
Covering an approximate square by \(O(1)\) standard squares changes the
corresponding counting estimate only by a fixed constant. For the tube
fibres, partition the fixed compact line-parameter window into finitely
many normal-coordinate charts. On each chart, the metric
\(d_{\mathcal A}\) and the Euclidean metric in the tube coordinates
used in \cite{DemeterORegan} are uniformly bi-Lipschitz equivalent.
Consequently, line-parameter balls are mapped into balls of comparable
radius, and the \((\delta,\sigma_s)\)-Katz--Tao condition is preserved,
up to a fixed change in its constant. Likewise, a family separated at
scale \(c\bar\delta\), for a fixed \(c>0\), can be partitioned into
\(O_c(1)\) \(\bar\delta\)-separated subfamilies by a bounded-degree
graph-colouring argument in the two-dimensional line-parameter space.
Fixed changes in tube width are absorbed by the choice of
\(\bar\delta\) and the preceding bounded coverings.

We may therefore apply \cite[Theorem~1.5]{DemeterORegan} separately in
each line chart and to each separation class. Since the number of
classes is fixed, H\"older's inequality gives
\[
\sum_{\nu} |\mathcal T_\nu|^{2/3}
\leq
O(1)\,|\mathcal T|^{2/3}.
\]
Finally, \(\bar\delta\asymp\delta\), and the notation
\(\lessapprox\) in \cite{DemeterORegan} means that, for every
\(\eta>0\), one may insert a factor \(\bar\delta^{-\eta}\).
Absorbing all fixed multiplicities and comparison constants therefore
gives
\[I(\mathcal Q,\mathcal T)
\lesssim_{s,\eta,C_Q,C_Y}
\delta^{-\eta-2s/3}
|\mathcal T|^{2/3}|\mathcal Q|^{1/3},\]
with additional dependence only on the fixed physical and
line-parameter windows.
\subsection{Rectangular non-concentration}

For $y\in M_2$, set
\[
 A_y:=\{F(x,y):x\in M_1\}.
\]

\begin{lemma}\label{lem:inc-row-rectangle}
For every interval $J$ with $|J|\geq\delta$ and every $y\in M_2$,
\begin{equation}
 |A_y\cap J|
 \ll\delta^{-\kappa\alpha}
 \left(\frac{|J|}{\delta}\right)^\alpha.
 \label{eq:inc-row-Frostman}
\end{equation}
Moreover, if $B$ is a rectangle of arbitrary orientation and side
lengths $r\geq r'\geq\delta$, then
\begin{equation}
 \#(P\cap B)
 \ll\delta^{-\kappa\alpha}
 \left(\frac{\sqrt{rr'}}{\delta}\right)^{2\alpha}.
 \label{eq:inc-P-rectangular}
\end{equation}
The constants depend only on the explicitly given parameters in Theorem~\ref{our_incidence_theorem}.
\end{lemma}

\begin{proof}
Fix $y\in M_2$ and an interval $J$, and let
\[
 E:=\{x\in M_1:F(x,y)\in J\}.
\]
The lower bound in \eqref{eq:inc-Phi-bilip} implies that
$x\mapsto F(x,y)$ is injective on $M_1$.  Hence
$|A_y\cap J|=|E|$.  If $E=\varnothing$, then the estimate \eqref{eq:inc-row-Frostman} is immediate. If $E\neq\varnothing$, choose $x_0\in E$.  For
every $x\in E$, the same lower bound in \eqref{eq:inc-Phi-bilip} gives
\[
 |x-x_0|
 \leq C_0\delta^{-\kappa}|F(x,y)-F(x_0,y)|
 \leq C_0\delta^{-\kappa}|J|.
\]
Thus $E$ lies in an interval of length at most
$2C_0\delta^{-\kappa}|J|$.  Applying
\eqref{eq:inc-M-Frostman} to $M_1$, we obtain
\[|E|=|M_1\cap E|\leq C_0{\left(\frac{2C_0\delta^{-\kappa}|J|}{\delta}\right)}^\alpha=C\delta^{-\kappa\alpha}{\left(\frac{|J|}{\delta}\right)}^\alpha,\]
where the constant $C=2^\alpha C_0^{\alpha+1}$. This proves
\eqref{eq:inc-row-Frostman}.

Let $H_B$ be the length of the vertical projection of $B$, and let
$W_B$ be the maximal length of a horizontal section of $B$.  If the
long side of $B$ makes an angle $\vartheta\in[0,\pi/2]$ with the
horizontal axis, then
\[ H_B=r\sin\vartheta+r'\cos\vartheta,\qquad W_B=\min\left\{\frac r{\cos\vartheta},\frac {r'}{\sin\vartheta}\right\},\]
with the usual endpoint interpretation.  Consequently, 
\begin{equation}
 H_B,W_B\geq r'\geq\delta,
 \qquad H_BW_B\leq2rr'.
 \label{eq:inc-HW}
\end{equation}

By \eqref{eq:inc-M-Frostman}, the vertical projection of $B$ meets at
most $C_0(H_B/\delta)^\alpha$ elements of $M_2$.  For each such $y$, the
horizontal section of $B$ is contained in an interval of length at most $W_B$; hence \eqref{eq:inc-row-Frostman} bounds its intersection with the row
$A_y\times\{y\}$ by
$C\delta^{-\kappa\alpha}(W_B/\delta)^\alpha$.  Multiplying these
bounds and using \eqref{eq:inc-HW} proves
\eqref{eq:inc-P-rectangular}.
\end{proof}

\subsection{Two regularization lemmas}

The factor $\delta^{-\kappa\alpha}$ in
\eqref{eq:inc-P-rectangular} is removed by decomposing the points into a
controlled number of rectangular Katz--Tao layers.

\begin{lemma}\label{lem:inc-rectangular-layers}
Fix $B_0,C_1\geq1$ and $\tau>0$.  Assign a nonnegative integer weight
$w(q)$ to each $\delta$-square in a fixed bounded square.  Suppose
$1\leq H\leq\delta^{-B_0}$, $w(q)\leq H$, and
\begin{equation}
 \sum_{c(q)\in B}w(q)
 \leq C_1H\left(\frac{\sqrt{rr'}}{\delta}\right)^{2\alpha}
 \label{eq:inc-weighted-rectangular}
\end{equation}
for every rectangle $B$ of arbitrary orientation and side lengths
$r\geq r'\geq\delta$.  Then, for all sufficiently small $\delta$, the
weighted points can be partitioned into
\begin{equation}
 J\leq C_{\alpha,\tau,B_0,C_1}H\delta^{-\tau}
 \label{eq:inc-point-layers-number}
\end{equation}
layers such that each layer contains at most one point in each
$\delta$-square and its occupied squares form a rectangular
$(\delta,2\alpha,C_{\alpha,\tau,B_0,C_1})$-Katz--Tao collection.
\end{lemma}

\begin{proof}
Take
\[J:=\left\lceil C_1H\delta^{-\tau}\right\rceil.\]
For every square $q$, regard the $w(q)$ points inside $q$ as labelled, and assign these points to a uniformly random injection into $\{1,\ldots,J\}$.  Make these choices independently for distinct squares.  Fix a colour $j$ and a rectangle $B$. Let 
\[X_{j,q}:=\begin{cases}
1&\text{if there is a point in }q\text{ receiving the colour }j,\\
0&\text{if there is no point in }q\text{ receiving the colour }j.
\end{cases}\]
If
$X_{j,B}$ is the number of squares with centre in $B$ that receive
colour $j$, then
\[
 X_{j,B}=\sum_{c(q)\in B}X_{j,q},
 \qquad
 \mathbb P(X_{j,q}=1)=\frac{w(q)}J,
\]
where the random variables $X_{j,q}$ are independent.  It follows from
\eqref{eq:inc-weighted-rectangular} that
\begin{equation}
 \mu_{j,B}:=\mathbb E X_{j,B}=\sum_{c(q)\in B}\mathbb{E}X_{j,q}
 \leq J^{-1}C_1H\left(\frac{\sqrt{rr'}}{\delta}\right)^{2\alpha}\leq\delta^\tau
 \left(\frac{\sqrt{rr'}}{\delta}\right)^{2\alpha}.
 \label{eq:inc-point-layer-mean}
\end{equation}
We now reduce the verification of the rectangular Katz--Tao condition to a finite family of test rectangles. Let
\[\mathcal{Q}_+ :=\{q : w(q)>0\}.\]
By hypothesis, all squares under consideration lie in a fixed bounded square. After replacing this square by a fixed enlargement if necessary, we may choose a fixed bounded square $S$ such that
\[c(q)\in S\quad \text{ for all }\quad q\in \mathcal{Q}_+.\]
We claim that there is a family \(\mathfrak{R}\) of rectangles satisfying
\[|\mathfrak{R}|=O\bigl(\delta^{-3}\log^2(1/\delta)\bigr)\]
with the following property: for every admissible rectangle \(B\) with $B\cap S\neq\varnothing$ of arbitrary orientation and side lengths \(r\ge r'\ge \delta\), there exists \(B^*\in\mathfrak R\) such that
\begin{equation}\label{eq:finite_cover}
B\cap\{c(q):q\in\mathcal{Q}_+\}\subset B^*,\qquad|B^*|\ll|B|.
\end{equation}
where $|\cdot|$ denotes Euclidean area. Note that if $B\cap S=\varnothing$, then the desired counting estimate is trivial. Thus, we focus on the case $B\cap S\neq\varnothing$ here.

Replace \(B\) by the smallest rectangle with the same orientation, side lengths at least \(\delta\), and containing \(B\cap S\). The resulting rectangle has centre in a fixed enlargement of \(S\) and area bounded above by \(|B|\). We then round its two side lengths upward to dyadic values, round its centre to a \(\delta\)-net in the fixed enlargement of \(S\), and round its angle to a \(\delta\)-net in \([0,\pi)\). For $\delta$ small enough, the resulting rectangle $B^*$ will satisfy \eqref{eq:finite_cover}. Let $\mathfrak{R}$ be the collection of all possible resulting rectangles $B^*$. Enlarging the rounded rectangle by a fixed multiplicative constant absorbs these perturbations, since the side length is at least \(\delta\). The number of possible centres is \(O(\delta^{-2})\), the number of possible angles is \(O(\delta^{-1})\), and the number of dyadic choices for the two side lengths is \(O(\log^2(1/\delta))\) since both the side lengths are $O_S(1)$. This proves the asserted cardinality bound $|\mathfrak{R}|=O\bigl(\delta^{-3}\log^2(1/\delta)\bigr)$ and the covering property \eqref{eq:finite_cover}.

For $B\in\mathfrak R$, set
\[A_B:=\left(\frac{\sqrt{rr'}}{\delta}\right)^{2\alpha}\geq1,\]
where $r,r'$ are the side lengths of $B$. Since $X_{j,B}$ is a Poisson-binomial random variable, it satisfies the
tail bound
\[
 \mathbb P(X_{j,B}\geq u)
 \leq\left(\frac{e\mu_{j,B}}{u}\right)^u
 \qquad\text{whenever }u\geq\mu_{j,B}.
\]
Together with \eqref{eq:inc-point-layer-mean}, this implies that, for a
sufficiently large fixed $C_2$ and all sufficiently small $\delta$,
\begin{equation}
 \mathbb P\bigl(X_{j,B}>C_2A_B\bigr)
 \leq\delta^{(\tau C_2/2)A_B}
 \leq\delta^{\tau C_2/2}.
 \label{eq:inc-point-layer-tail}
\end{equation}
Let $E_{j,B}$ denote the event $X_{j,B}>C_2A_B$ and define the event 
\[\mathcal{E}:=\bigcup_{\substack{1\leq j\leq J\\B\in\mathfrak{R}}}E_{j,B}.\]
The number of possible pairs $(j,B)$ above is at most $J|\mathfrak{R}|\ll\delta^{-A_1}$ for some $A_1=A_1(B_0,\tau)$. Therefore, 
\[\mathbb{P}(\mathcal{E})\leq\delta^{-A_1}\max_{j,B}\mathbb{P}(E_{j,B})\leq\delta^{\frac{\tau C_2}{2}-A_1}.\]
Thus, we may choose $C_2>\frac{2(A_1+1)}{\tau}$ so that $\mathbb{P}(\mathcal{E})\leq\delta<1$. That is, there is a deterministic colouring satisfying $X_{j,B}\leq C_2A_B$ for every colour $1\leq j\leq J$ and $B\in\mathfrak{R}$. The injection used in each square gives the
first conclusion of the lemma. 

Next, via \eqref{eq:finite_cover}, we know that for every rectangle $B$ and colour $j$, there exists a canonical rectangle $B^*\in\mathfrak{R}$ such that 
\[\#\{q:c(q)\in B,q\text{ receives the colour }j\}\leq X_{j,B^*}\leq C_2A_{B^*}\ll C_2A_B\]
This transfers the bound from the canonical rectangles to all admissible rectangles and completes the proof.
\end{proof}

The second lemma regularizes, separately at each degree, the family of
lines through a square.

\begin{lemma}
\label{lem:inc-degree-colouring}
Fix $B_0\geq1$, $0<\sigma\leq1$, $\tau>0$, a bounded square $S\subset\mathbb{R}^2$, and a compact line window $\Omega\subset\mathcal{A}$.  Let $0<\delta<1$ and $\mathcal L$ be a
$\delta$-separated family of $\delta$-tubes whose axial lines lie in $\Omega$.  Let $D\geq1$ and $\mathcal Q_D$ be a
collection of $\delta$-squares $q$ contained in $S$ and
satisfying
\begin{equation}
 D\leq\deg_{\mathcal L}(q)<2D,\qquad\deg_{\mathcal L}(q):=\#\{T\in\mathcal L:q\cap T\neq\varnothing\}.
 \label{eq:inc-degree-class}
\end{equation}
Suppose that
$|\mathcal L|+|\mathcal Q_D|\leq\delta^{-B_0}$.  Then, for all
sufficiently small $\delta$, the tube family $\mathcal L$ can be partitioned into
\begin{equation}
 K\leq C_{\tau,B_0,\sigma,S,\Omega}\delta^{-\tau}D^{1-\sigma}
 \label{eq:inc-line-colours-number}
\end{equation}
classes $\mathcal L_k$ such that, for every $q\in\mathcal Q_D$, the
axial lines of the incident family
\[
 \mathcal L_k(q):=\{T\in\mathcal L_k:q\cap T\neq\varnothing\}
\]
form a $(\delta,\sigma,C_{\tau,B_0,\sigma,S,\Omega})$-Katz--Tao family in
$\mathcal A$.
\end{lemma}

\begin{proof}
Let $\mathcal{L}(q):=\{T\in\mathcal{L}:q\cap T\neq\emptyset\}$. We first prove that for every $q\in\mathcal{Q}_D,u\in\mathcal{A}$, and $r\geq\delta$, 
\[\#(\mathcal{L}(q)\cap B_\mathcal{A}(u,r))\lesssim_{S,\Omega}\min\{D,r/\delta\}.\]
Let ${\{(U_\nu,\psi_\nu)\}}_{\nu}$ be the charts of $\mathcal{A}$. Cover the compact line window $\Omega$ by finitely many relatively compact normal-coordinate charts $V_\nu\Subset U_\nu$. On each \(U_\nu\), choose a continuous unit normal vector \(\omega_\nu(\theta)\) for the lines in the chart, so that every \(\ell\in U_\nu\) can be written uniquely as \(\ell=\ell_{\omega_\nu(\theta),b}\) and $\psi_\nu(\ell)=(\theta,b)$ in the corresponding local coordinates \((\theta,b)\). Additionally, the charts $V_\nu$ may be chosen so that there is an absolute constant $C\geq1$ satisfying 
\[C^{-1}|\psi_\nu(\ell)-\psi_\nu(\ell')|\leq d_\mathcal{A}(\ell,\ell')\leq C|\psi_\nu(\ell)-\psi_\nu(\ell')|\quad\text{ whenever}\quad\ell,\ell'\in V_\nu.\]
We may do a refinement to the original charts and make the charts ${\{(U_\nu,\psi_\nu)\}}_\nu$ disjoint. Note that the number of charts is still a constant depending on $\Omega$ after this refinement.

Let $c(q)$ be the centre of $q$. If the $\delta$-tube about $\ell_{\omega,b}$ meets $q$, choose $z\in q\cap\ell_{\omega,b}^{(\delta)}$. Since $|z-c(q)|\leq\sqrt{2}\delta/2$,
\[|b-c(q)\cdot\omega|\leq|b-z\cdot\omega|+|(z-c(q))\cdot\omega|\leq\left(1+\frac{\sqrt{2}}{2}\right)\delta.\]
Thus, in the $\nu$-th chart, the parameters of lines incident to $q$ lie in a $C\delta$-neighbourhood of the graph
\[
 \Gamma_{q,\nu}:=\{(\theta,g_{q,\nu}(\theta)):g_{q,\nu}(\theta)=c(q)\cdot\omega_\nu(\theta)\}.
\]
Since $c(q)$ ranges over the fixed bounded square $S$ and there are only finitely many charts, the functions $g_{q,\nu}$ have Lipschitz constants bounded uniformly in $q$ and $\nu$.

Fix $u\in\mathcal A$, $r\geq\delta$, and one of the charts $V_\nu$. If
\[
 \mathcal L(q)\cap B_{\mathcal A}(u,r)\cap V_\nu\neq\varnothing,
\]
choose one of its axial lines $\ell_0$. Every other axial line $\ell$ in this set satisfies $d_{\mathcal A}(\ell,\ell_0)\leq2r$. By the coordinate comparison above, the corresponding parameters lie in an $O_\Omega(r)$-ball about $\psi_\nu(\ell_0)$ and, by the preceding incidence estimate, in a $C\delta$-neighbourhood of the uniformly Lipschitz graph $\Gamma_{q,\nu}$. This portion of the graph neighbourhood can therefore be covered by
\[
 O_{S,\Omega}(1+r/\delta)=O_{S,\Omega}(r/\delta)
\]
coordinate balls of radius $C'\delta$. Each such coordinate ball contains only $O_\Omega(1)$ axial lines of $\mathcal L$: by the bi-Lipschitz comparison of the chart metric with $d_{\mathcal A}$, it is contained in a $d_{\mathcal A}$-ball of radius $C''\delta$, and a $\delta$-separated subset of a fixed two-dimensional chart has $O_\Omega(1)$ points in such a ball. Summing over the finitely many charts gives
\[
 \#\bigl(\mathcal L(q)\cap B_{\mathcal A}(u,r)\bigr)
 \leq C_{\mathrm{pack}}\frac r\delta,
\]
where $C_{\mathrm{pack}}$ is a constant depending only on $S$ and $\Omega$. On the other hand, \eqref{eq:inc-degree-class} gives $\#\mathcal L(q)<2D$. After enlarging $C_{\mathrm{pack}}$ if necessary, we obtain
\begin{equation}
 \#\bigl(\mathcal L(q)\cap B_{\mathcal A}(u,r)\bigr)
 \leq C_{\mathrm{pack}}\min\{D,r/\delta\}.
 \label{eq:inc-raw-fibre-packing}
\end{equation}

We next colour the tubes. Set
\[
 \Lambda:=\delta^{-\tau}D^{1-\sigma},\qquad K_0:=\lceil\Lambda\rceil.
\]
Since $\Lambda\geq1$,
\begin{equation}
 K_0\leq2\delta^{-\tau}D^{1-\sigma}.
 \label{eq:inc-K0-bound}
\end{equation}
 If $K_0>|\mathcal L|$, put each tube in a separate class. The number of classes is then smaller than $K_0$, and every incident colour fibre contains at most one tube, so
\[
 \#\bigl(\mathcal L_k(q)\cap B_{\mathcal A}(u,r)\bigr)
 \leq1\leq(r/\delta)^\sigma
 \qquad(r\geq\delta).
\]
Thus the conclusion holds in this case. We may therefore assume $K_0\leq|\mathcal L|$, put $K=K_0$, and colour the tubes independently and uniformly with colours $1,\ldots,K$.

Fix $q\in\mathcal Q_D$, $u\in\mathcal A$, $r\geq\delta$, and a colour $k$, and write $x=r/\delta\geq1$. Define
\[
 X(q,u,r,k):=\#\bigl(\mathcal L_k(q)\cap B_{\mathcal A}(u,r)\bigr).
\]
This is a binomial random variable. By \eqref{eq:inc-raw-fibre-packing} and $K\geq\Lambda$, its mean satisfies
\begin{equation}
 \mu(q,u,r,k)
 \leq \frac{C_{\mathrm{pack}}\min\{D,x\}}{K}
 \leq C_{\mathrm{pack}}\delta^\tau\min\{D,x\}D^{\sigma-1}
 \leq C_{\mathrm{pack}}\delta^\tau x^\sigma.
 \label{eq:inc-line-colour-mean}
\end{equation}
For the last inequality, if $x\leq D$, then $xD^{\sigma-1}\leq x^\sigma$, while if $x\geq D$, then $D^\sigma\leq x^\sigma$.

It remains to obtain this estimate simultaneously for all centres and radii. Put
\[
 R_\Omega:=\max\{1,\operatorname{diam}_{\mathcal A}(\Omega)\},\qquad r_j:=2^j\delta,
\]
and let $J_\Omega$ be the smallest integer for which $r_{J_\Omega}\geq R_\Omega$. Then
\begin{equation}
 J_\Omega<\log_2(R_\Omega/\delta)+1.
 \label{eq:inc-dyadic-radii}
\end{equation}
For $T\in\mathcal L$, write $\ell_T$ for its axial line. Initially test only the balls
\[
 B_{\mathcal A}(\ell_T,r_j),
 \qquad T\in\mathcal L,\quad0\leq j\leq J_\Omega.
\]
There are at most
\[
 |\mathcal Q_D|\,|\mathcal L|\,K(J_\Omega+1)
\]
choices of $(q,T,k,j)$. In the random-colouring case $K\leq|\mathcal L|$, so the polynomial cardinality hypothesis and \eqref{eq:inc-dyadic-radii} imply that, for all sufficiently small $\delta$, the total number of tested events is at most
\[
 \delta^{-A},\qquad A:=3B_0+2.
\]
where $B_0$ is the bounding constant of both $|\mathcal{L}|$ and $|\mathcal{Q}_D|$ mentioned in the assumption.
Fix a constant $C_*>0$, to be chosen below. For a tested radius set $x_j:=r_j/\delta$ and
\[
 m_j:=\lfloor C_*x_j^\sigma\rfloor+1.
\]
By \eqref{eq:inc-line-colour-mean}, for all sufficiently small $\delta$, uniformly over all tested parameters,
\[
 m_j\geq\mu(q,\ell_T,r_j,k),
 \qquad
 \frac{e\mu(q,\ell_T,r_j,k)}{m_j}\leq\frac{eC_\text{pack}\delta^\tau x_j^\sigma}{C_*x_j^\sigma}\leq\delta^{\tau/2}.
\]
The Poisson--binomial tail estimate used above therefore gives
\begin{align*}
 \mathbb P\Bigl(X(q,\ell_T,r_j,k)>C_*x_j^\sigma\Bigr)
 &\leq
 \left(\frac{e\mu(q,\ell_T,r_j,k)}{m_j}\right)^{m_j}\\
 &\leq\delta^{\tau m_j/2}
 \leq\delta^{(\tau C_*/2)x_j^\sigma}
 \leq\delta^{\tau C_*/2}.
\end{align*}
Choose
\[
 C_*>\frac{2(A+1)}{\tau}.
\]
The union bound shows that the probability that any tested estimate fails is at most
\[
 \delta^{-A}\delta^{\tau C_*/2}<1.
\]
Hence there is a deterministic colouring for which
\begin{equation}
 \#\bigl(\mathcal L_k(q)\cap B_{\mathcal A}(\ell_T,r_j)\bigr)
 \leq C_*\left(\frac{r_j}{\delta}\right)^\sigma
 \label{eq:inc-tested-ball-colour}
\end{equation}
for every $q\in\mathcal Q_D$, $T\in\mathcal L$, $1\leq k\leq K$, and $0\leq j\leq J_\Omega$.

Finally, transfer the tested bounds to arbitrary balls. Fix $u\in\mathcal A$, $r\geq\delta$, $q\in\mathcal Q_D$, and a colour $k$. If $\mathcal L_k(q)\cap B_{\mathcal A}(u,r)=\varnothing$, the desired estimate is trivial. Otherwise choose an axial line $v$ in this intersection. The triangle inequality gives
\[
 \mathcal L_k(q)\cap B_{\mathcal A}(u,r)
 \subset \mathcal L_k(q)\cap B_{\mathcal A}(v,2r).
\]
If $2r\leq r_{J_\Omega}$, choose $j$ so that $2r\leq r_j<4r$. Since $v$ is the axial line of a member of $\mathcal L$, \eqref{eq:inc-tested-ball-colour} yields
\[
 \#\bigl(\mathcal L_k(q)\cap B_{\mathcal A}(u,r)\bigr)
 \leq C_*4^\sigma(r/\delta)^\sigma.
\]
If $2r>r_{J_\Omega}$, then $B_{\mathcal A}(v,r_{J_\Omega})$ contains all of $\Omega$, because $r_{J_\Omega}\geq\operatorname{diam}_{\mathcal A}(\Omega)$. Hence \eqref{eq:inc-tested-ball-colour}, applied at the largest tested radius, gives
\[
 \#\bigl(\mathcal L_k(q)\cap B_{\mathcal A}(u,r)\bigr)
 \leq C_*\left(\frac{r_{J_\Omega}}{\delta}\right)^\sigma
 < C_*2^\sigma(r/\delta)^\sigma.
\]
Since $0<\sigma\leq1$, both cases imply
\begin{equation}
 \#\bigl(\mathcal L_k(q)\cap B_{\mathcal A}(u,r)\bigr)
 \leq4C_*(r/\delta)^\sigma
 \qquad(u\in\mathcal A,\ r\geq\delta).
 \label{eq:inc-arbitrary-ball-transfer}
\end{equation}
Every $\mathcal L_k$ is a subfamily of the original $\delta$-separated family $\mathcal L$, hence remains $\delta$-separated. Thus \eqref{eq:inc-arbitrary-ball-transfer} is precisely the Katz--Tao condition \eqref{eq:inc-line-KT}, with a constant depending only on $\tau,B_0,\sigma,S$, and $\Omega$.

Finally, \eqref{eq:inc-K0-bound} gives
\[
 K\leq2\delta^{-\tau}D^{1-\sigma}
\]
in the random-colouring case, and the same bound holds in the singleton case. This proves \eqref{eq:inc-line-colours-number} and completes the proof.

\end{proof}

\begin{proof}[Proof of Theorem~\ref{our_incidence_theorem}]
If either $M_1$ or $M_2$ is empty, the conclusion is immediate.  The
lower bound in \eqref{eq:inc-Phi-bilip} shows that $\Phi$ is injective
on $M_1\times M_2$; in particular,
\begin{equation}
 |P|=|M_1||M_2|\leq C_0^2\delta^{-2\alpha}.
 \label{eq:inc-P-cardinality}
\end{equation}

Discard every line whose $\delta$-neighbourhood misses $P$.  By
Remark~\ref{rem:inc-bounded-window}, the remaining axial lines lie in a
fixed compact subset of $\mathcal A$.  Assign each point of $P$ to one
of the half-open squares of the standard $\delta$-grid that contains it,
and let $\mathcal L$ be the corresponding family of $\delta$-tubes.  If
$p\in\ell^{(\delta)}$, then the square assigned to $p$ meets the tube
about $\ell$.  Thus every original point-line incidence gives a
square-tube incidence, without changing the scale.

For a grid square $q$, define
\[
 w(q):=\#\{p\in P:p\text{ is assigned to }q\}.
\]
If $c(q)$ lies in a rectangle $B$ with side lengths at least $\delta$,
then $q$ is contained in an $O(\delta)$-enlargement of $B$.  That
enlargement can be covered by a fixed number of rectangles whose side
lengths are comparable to those of $B$.  Lemma~\ref{lem:inc-row-rectangle}
therefore gives a constant $C$ such that 
\begin{equation}
 w(q)\leq C\delta^{-\kappa\alpha}
 \label{eq:inc-square-weight}
\end{equation}
and
\begin{equation}
 \sum_{c(q)\in B}w(q)
 \leq C\delta^{-\kappa\alpha}
 \left(\frac{\sqrt{rr'}}{\delta}\right)^{2\alpha}.
 \label{eq:inc-weighted-P-bound}
\end{equation}

Apply Lemma~\ref{lem:inc-rectangular-layers} with
\[
 H=C\delta^{-\kappa\alpha}.
\]
After increasing $C$ and decreasing $\delta_0$, we may ensure $H$ satisfies $1\leq H\leq\delta^{-2}$.  Consequently, for every $\tau_p>0$, there exists $C_{\tau_p}=C_{\tau_p}(\alpha,C)$ such that the labelled points of $P$ can be partitioned into
\begin{equation}
 J\leq C_{\tau_p}\delta^{-\kappa\alpha-\tau_p}
 \label{eq:inc-J}
\end{equation}
layers.  Each layer contains at most one labelled point in a square, and its occupied squares form a rectangular $(\delta,2\alpha,C_{\tau_p})$-Katz--Tao collection.  Point-tube incidences in a layer are bounded by square-tube incidences, and summing over the layers counts all original incidences.

Fix one layer $\mathcal Q$.  Decompose its squares according to their
degree in the full tube family $\mathcal L$:
\[
 \mathcal Q_D:=\{q\in\mathcal Q:
                  D\leq\deg_{\mathcal L}(q)<2D\},
 \qquad n_D:=|\mathcal Q_D|,
\]
where $D$ ranges over dyadic integers.  Degree-zero squares contribute
nothing.  By \eqref{eq:inc-raw-fibre-packing}, $D\ll\delta^{-1}$.

Apply Lemma~\ref{lem:inc-degree-colouring} with $\sigma$ from
\eqref{eq:inc-sigma} and an arbitrary $\tau_\ell>0$.  Its polynomial
cardinality hypothesis holds with a fixed $B_0$: the bounded physical
window contains $O(\delta^{-2})$ grid squares and
$|\mathcal L|\leq C_0\delta^{-2\alpha}\leq C_0\delta^{-2}$.  Hence
$\mathcal L$ is the disjoint union of
\begin{equation}
 K\leq C_{\tau_\ell,B_0,\sigma}
 \delta^{-\tau_\ell}D^{1-\sigma}
 \label{eq:inc-K}
\end{equation}
colours $\mathcal L_k$, where the constant may also depend on the fixed
physical and line-parameter windows.  Each incident colour fibre has the
Katz--Tao property required in
Theorem~\ref{thm:inc-DO-rectangular}. Since a subcollection of a
rectangular Katz--Tao collection satisfies the same upper bound, Theorem \ref{thm:inc-DO-rectangular} applies to every pair $(\mathcal Q_D,\mathcal L_k)$.  For every
$\eta>0$, H\"older's inequality in $k$ gives
\begin{align*}
 Dn_D
 &\leq I(\mathcal Q_D,\mathcal L)
   =\sum_{k=1}^K I(\mathcal Q_D,\mathcal L_k) \notag\\
 &\leq C_{\alpha,\eta,\tau_p,\tau_\ell,B_0}
 \delta^{-\eta-2\alpha/3}n_D^{1/3}
 \sum_{k=1}^K|\mathcal L_k|^{2/3} \notag\\
 &\leq C_{\alpha,\eta,\tau_p,\tau_\ell,B_0}
 \delta^{-\eta-2\alpha/3}n_D^{1/3}
 K^{1/3}|\mathcal L|^{2/3} \notag\\
 &\leq \underbrace{C_{\alpha,\eta,\tau_p,\tau_\ell,B_0}C_{\tau_\ell,B_0,\sigma}^{1/3}C_0}_{=:C_*}
 \delta^{-\eta-2\alpha-\tau_\ell/3}
 D^{(1-\sigma)/3}n_D^{1/3}.
\end{align*}
Therefore, solving for $n_D$, we obtain
\begin{equation}\label{eq:inc-degree}
 n_D
 \leq C_*^{3/2}
 \delta^{-3\alpha-3\eta/2-\tau_\ell/2}
 D^{-1-\sigma/2}.
\end{equation}
On the other hand,
\[I(\mathcal Q_D,\mathcal L)<2Dn_D.\]
Together with \eqref{eq:inc-degree} for $n_D$, we get
\begin{equation}
 I(\mathcal Q_D,\mathcal L)
 \leq 2C_*^{3/2}\delta^{-3\alpha-3\eta/2-\tau_\ell/2}D^{-\sigma/2}.
 \label{eq:inc-degree-solved}
\end{equation}

We now separate low and high degrees.  Define
\begin{equation}
 \zeta:=\frac{G}{2\sigma}
 =\frac\alpha4-\frac{\alpha\kappa}{2\sigma},
 \qquad
 D_0:=\delta^{-\alpha+\zeta}.
 \label{eq:inc-D0}
\end{equation}
Since $0\leq\kappa<\sigma/2$,
\begin{equation}
 0<\zeta\leq\frac\alpha4,
 \qquad
 \alpha-\zeta\geq\frac{3\alpha}{4}>0.
 \label{eq:inc-zeta-range}
\end{equation}
Choose a dyadic integer $D_*$ with $D_0\leq D_*<2D_0$.  The points
assigned to squares of degree smaller than $D_*$ contribute, by
\eqref{eq:inc-P-cardinality}, at most
\begin{equation}
 I_{\mathrm{low}}
 \leq D_*|P|
 \leq C_0^2\delta^{-3\alpha+\zeta}.
 \label{eq:inc-low-degree}
\end{equation}

For $D\geq D_*$, sum \eqref{eq:inc-degree-solved} over the $J$ point
layers and over the $O(\log(1/\delta))$ dyadic degree classes.  Given any $\rho>0$, the logarithmic factor is at most $C_\rho\delta^{-\rho}$.  Equations \eqref{eq:inc-J},
\eqref{eq:inc-D0}, and \eqref{eq:inc-degree-solved} give
\begin{align}
 I_{\mathrm{high}}
 &\leq
 C_\rho\delta^{-\rho}J
 \max_{\substack{D\geq D_*\\ D\ \mathrm{dyadic}}}
 I(\mathcal Q_D,\mathcal L)
 \notag\\
 &\leq
 2C_\rho C_{\tau_p}C_*^{3/2}
 \delta^{-3\alpha-3\eta/2-\tau_\ell/2-\rho-\kappa\alpha-\tau_p}
 D_*^{-\sigma/2}
 \notag\\
 &\leq
 2C_\rho C_{\tau_p}C_*^{3/2}
 \delta^{-3\alpha+G-\zeta\sigma/2-\tau_p
 -3\eta/2-\tau_\ell/2-\rho}.
 \label{eq:inc-high-degree}
\end{align}
Here we used the identity
$G=\alpha\sigma/2-\alpha\kappa$.  Since
$\zeta\sigma/2=G/4$, choose
$\tau_p,\tau_\ell,\eta,\rho>0$ so small that
\begin{equation}
 \tau_p+\frac{3\eta}{2}+\frac{\tau_\ell}{2}+\rho<\frac G4.
 \label{eq:inc-loss-choice}
\end{equation}
Since $\zeta\sigma/2=G/4$, the exponent in
\eqref{eq:inc-high-degree} equals
\[
 -3\alpha+\frac{3G}{4}
 -\left(
   \tau_p+\frac{3\eta}{2}+\frac{\tau_\ell}{2}+\rho
  \right).
\]
Thus \eqref{eq:inc-loss-choice} gives
\begin{equation}
 I_{\mathrm{high}}
 \leq
 2C_\rho C_{\tau_p}C_*^{3/2}
 \delta^{-3\alpha+G/2}.
 \label{eq:inc-high-final}
\end{equation}
Moreover, since $\sigma\leq1$,
\[
 \zeta=\frac{G}{2\sigma}\geq\frac{G}{2}.
\]
Hence \eqref{eq:inc-low-degree} also gives
\[
 I_{\mathrm{low}}\lesssim\delta^{-3\alpha+G/2}.
\]
Combining this with \eqref{eq:inc-high-final}, and using
$0<\varepsilon<G/4<G/2$ and $0<\delta<1$, we obtain
\[
 I_{\mathrm{low}}+I_{\mathrm{high}}
 \lesssim
 \delta^{-3\alpha+G/2}
 \leq
 \delta^{-3\alpha+\varepsilon}.
\]
This proves \eqref{eq:inc-main-bound}.
\end{proof}

We shall use the theorem with fixed changes in the incidence width and
the separation scale.  The following corollary makes this stability
explicit.

\begin{corollary}
\label{cor:inc-fixed-scale-robustness}
Fix $c,C_1>0$.  Under all the hypotheses of
Theorem~\ref{our_incidence_theorem}, except that the line family $L$ is
only assumed to be $c\delta$-separated, one has
\begin{equation}
 \#\{(p,\ell)\in P\times L:p\in\ell^{(C_1\delta)}\}
 \leq C'\delta^{-3\alpha+\varepsilon}.
 \label{eq:inc-robust-bound}
\end{equation}
Here $C'$ and the upper endpoint for $\delta$ may also depend on $c$ and
$C_1$.  The conclusion remains valid if either $M_i$ is replaced by any
subset satisfying the same upper bounds in
\eqref{eq:inc-M-Frostman}.
\end{corollary}

\begin{proof}
The final assertion is already contained in Theorem \ref{our_incidence_theorem}, so it remains to prove \eqref{eq:inc-robust-bound}.  Choose a fixed integer $A\geq\max\{1,C_1\}$ and put $\delta':=A\delta$.  Discard every line whose $C_1\delta$-neighbourhood does not meet $P$. Since $P\subset B(0,R)$, every remaining line has distance at most $R+C_1\delta$ from the origin.  Consequently, for $\delta$ sufficiently small, all remaining lines lie in a fixed compact subset of $\mathcal A$ depending only on $R$ and $C_1$. Form a graph with vertex set $L$ by joining two lines whenever their $d_{\mathcal A}$-distance is less than $\delta'$.  Since $L$ is $c\delta$-separated and $\delta'=A\delta$, a packing argument in the two-dimensional line-parameter space shows that every $d_{\mathcal A}$-ball of radius $\delta'$ contains at most $O_{A,c}(1)$ elements of $L$.  The graph therefore has degree $O_{A,c}(1)$.  A greedy colouring partitions $L$ into $O_{A,c}(1)$ subfamilies, each of which is $\delta'$-separated.

After ordering the points of \(M_i\), partition them according to their indices modulo \(A\). Then, we obtain $O_A(1)$ subsets such that any two points in the same subset are separated by at least \(A\delta=\delta'\) since $M_i$ is $\delta$-separated. Consequently,
\[
 P=\bigcup_{a=1}^{A}\bigcup_{b=1}^{A}\Phi(M_{1,a}\times M_{2,b}).
\]
This union need not be disjoint, but summing the corresponding incidence
bounds over its at most $A^2$ product pieces is sufficient for an upper
bound.
For intervals $I$ with $|I|\geq\delta'$, these subsets satisfy
\[|M_{i,a}\cap I|
 \leq|M_i\cap I|\leq C_0{\left(\frac{|I|}{\delta}\right)}^\alpha=C_0A^\alpha\left(\frac{|I|}{\delta'}\right)^\alpha,
 \qquad
 |M_{i,a}|\leq|M_i|\leq C_0\delta^{-\alpha}=C_0A^\alpha(\delta')^{-\alpha}.\]
Furthermore, \eqref{eq:inc-Phi-bilip} becomes
\[
 (C_0A^\kappa)^{-1}(\delta')^\kappa|u-v|
 \leq|\Phi(u)-\Phi(v)|
 \leq C_0A^\kappa(\delta')^{-\kappa}|u-v|,
\]
and \eqref{eq:inc-line-cardinality} becomes
$|L|\leq C_0A^{2\alpha}(\delta')^{-2\alpha}$.  Finally, since $C_1\delta\leq\delta'$,
\[
 \ell^{(C_1\delta)}\subset\ell^{(\delta')}
 \qquad(\ell\in L).
\]
After decreasing the upper endpoint for $\delta$ so that $\delta'$ lies
in the admissible scale range of Theorem~\ref{our_incidence_theorem},
we may apply that theorem at scale $\delta'$, with the enlarged fixed
constant
\[
 C_0':=C_0\max\{A^\alpha,A^{2\alpha},A^\kappa\},
\]
therefore applies to every product piece and every line-colour class.
Summing over the at most $A^2$ product pieces and the
$O_{A,c}(1)$ line-colour classes gives
\[
 \#\{(p,\ell)\in P\times L:
       p\in\ell^{(C_1\delta)}\}
 \leq
 C_{A,c}(\delta')^{-3\alpha+\varepsilon}
 \leq
 C'\delta^{-3\alpha+\varepsilon},
\]
as required.
\end{proof}

\section{Proof of the main result (Theorem~\ref{thm-main-1})}

For a probability measure $\mu$ on $[0,1]$, write
$\mu_3:=\mu^{\otimes3}$ and define the coincidence function
\begin{equation}\label{eq:coincidence-function}
 \mathcal C_\mu(\rho)
 :=
 (\mu_3\times\mu_3)
 \Bigl(
  \bigl\{(u,v)\in[0,1]^3\times[0,1]^3:
  |f(u)-f(v)|\leq\rho\bigr\}
 \Bigr).
\end{equation}

\begin{lemma}\label{lem:schwartz-reduction}
Let $\eta>0$.  Suppose that there are $\rho_0>0$ and $C_1<\infty$ such
that
\begin{equation}\label{eq:coincidence-target}
 \mathcal C_\mu(\rho)\leq C_1\rho^{\alpha+\eta}
 \quad\text{ for all }\quad0<\rho\leq\rho_0.
\end{equation}
Then, for every fixed nonnegative Schwartz function $\varphi$,
\begin{equation}\label{eq:schwartz-reduction-conclusion}
 \|\varphi_\delta*f_\#\mu_3\|_{L^2(\mathbb R)}^2
 \lesssim_{\varphi,C_1,\rho_0}
 \delta^{\alpha+\eta-1}\quad\text{ for all }\quad0<\delta\leq\rho_0/2.
\end{equation}
\end{lemma}

\begin{proof}
Set $\widetilde\varphi(t):=\varphi(-t)$ and
$K:=\varphi*\widetilde\varphi$.  The function $K$ is a nonnegative
Schwartz function.  For every integer $N>\alpha+\eta$,
\begin{equation}\label{eq:schwartz-dyadic-majorant}
 K(s)\lesssim_{\varphi,N}
 \sum_{j=0}^{\infty}2^{-jN}\mathbf 1_{\{|s|\leq2^j\}}.
\end{equation}
If $\nu=f_\#\mu_3$, Tonelli's theorem gives the exact identity
\[
 \|\varphi_\delta*\nu\|_2^2
 =\iint \delta^{-1}K\left(\frac{f(u)-f(v)}{\delta}\right)
 \,d\mu_3(u)\,d\mu_3(v).
\]
Consequently, \eqref{eq:schwartz-dyadic-majorant} implies
\begin{equation}\label{eq:schwartz-dyadic-sum}
 \|\varphi_\delta*\nu\|_2^2
 \lesssim_{\varphi,N}
 \delta^{-1}\sum_{j=0}^{\infty}2^{-jN}
 \mathcal C_\mu(2^j\delta).
\end{equation}
For the indices satisfying $2^j\delta\leq\rho_0$, use
\eqref{eq:coincidence-target}.  Their contribution is
\[
 \lesssim C_1\delta^{\alpha+\eta-1}
 \sum_{j\geq0}2^{-j(N-\alpha-\eta)}
 \lesssim C_1\delta^{\alpha+\eta-1},
\]
provided $N>\alpha+\eta$. 

For the remaining indices, define
\[J_0:=\max\{j\geq0:2^j\delta\leq\rho_0\}.\]
This integer is well-defined because $0<\delta\leq\rho_0/2$.
For $j>J_0$, use the trivial bound
$\mathcal C_\mu(2^j\delta)\leq1$, which follows because $\mu_3$ is a
probability measure. Summing the resulting geometric series yields
\[\delta^{-1}\sum_{j>J_0}2^{-jN}\mathcal{C}_\mu(2^j\delta)\leq\delta^{-1}\sum_{j>J_0}2^{-jN}\lesssim_N\delta^{-1}2^{-(J_0+1)N}.\]
By the maximality of \(J_0\),
\[\rho_0<2^{J_0+1}\delta,\]
and hence
\[2^{-(J_0+1)N}\leq\left(\frac{\delta}{\rho_0}\right)^N.\]
Therefore, these remaining indices contribute
\[
\delta^{-1}\sum_{j>J_0}2^{-jN}\mathcal{C}_\mu(2^j\delta)=O_{\rho_0,N}(\delta^{N-1}).
\]
Choosing \(N>\alpha+\eta\), and using
\(0<\delta\leq\rho_0/2\), we have
\[
\delta^{N-1}
\lesssim_{\rho_0,N,\alpha,\eta}
\delta^{\alpha+\eta-1}.
\]
Combining both cases proves \eqref{eq:schwartz-reduction-conclusion}.
\end{proof}

\subsection{Proof of Theorem~\ref{thm-main-1} for
\texorpdfstring{$\alpha$-AD-regular}{alpha-AD-regular}
measures}

\begin{definition}\label{def:AD_alphab}
Let $\alpha\in(0,1)$.
We say that a probability measure $\mu$ is $\alpha$-AD-regular
on its support with constant $C_\mu>0$ if
$\operatorname{diam}(\spt\mu)>0$ and, for every
$x\in\spt(\mu)$ and every
$0<r\leq\operatorname{diam}(\spt\mu)$, one has
\[
 C_\mu^{-1}r^\alpha
 \leq
 \mu\bigl(B(x,r)\bigr)
 \leq
 C_\mu r^\alpha.
\]
\end{definition}

When $C_\mu$ is an absolute constant and the support of $\mu$ is clear
from the context, we simply say that $\mu$ is $\alpha$-AD-regular.

We first prove the coincidence estimate
\eqref{eq:coincidence-target} under the additional assumption that $\mu$
is $\alpha$-AD-regular.  The value of the exponent $\eta$ will be fixed
after the transverse and small-gradient bounds have been combined.

Put
\[
 \sigma:=\min\{2\alpha,2-2\alpha\}.
\]
Choose a small positive constant
\[
 \kappa<\frac{\alpha^2\sigma}{64},
\]
put
\[
 \rho:=\delta^\kappa,
\]
and define
\[
 \mathcal G_s
 :=
 \bigl\{
 u\in\operatorname{spt}(\mu)^3:
 |\nabla f(u)|\leq s
 \bigr\}.
\]
For $u=(x,y,z)$ and $v=(x',y',z')$, let
\[
\begin{array}{lll}
 D_1(u,v)=\partial_{x'}f(v),&
 D_2(u,v)=\partial_{y'}f(v),&
 D_3(u,v)=\partial_{z'}f(v),\\
 D_4(u,v)=\partial_xf(u),&
 D_5(u,v)=\partial_yf(u),&
 D_6(u,v)=\partial_zf(u).
\end{array}
\]
Set
\begin{align}
 I_0
 &:=
 \iint
 \mathbf 1_{\{|f(u)-f(v)|\leq2\delta\}}
 \mathbf 1_{\mathcal G_{2\rho}}(u)
 \mathbf 1_{\mathcal G_{2\rho}}(v)
 \,d\mu_3(u)\,d\mu_3(v),
 \label{eq:def_I0}\\
 I_j
 &:=
 \iint
 \mathbf 1_{\{|f(u)-f(v)|\leq2\delta\}}
 \mathbf 1_{\{|D_j(u,v)|>\rho\}}
 \,d\mu_3(u)\,d\mu_3(v),
 \qquad 1\leq j\leq6.
 \label{eq:def_Ij}
\end{align}
If
\[
 (u,v)\notin\mathcal G_{2\rho}\times\mathcal G_{2\rho},
\]
then one of the six quantities $|D_j(u,v)|$ is larger than $\rho$,
since $2/\sqrt3>1$.  Therefore,
\begin{equation}\label{eq:coincidence-six-term-decomposition}
 \mathcal C_\mu(2\delta)
 \leq
 I_0+\sum_{j=1}^6I_j.
\end{equation}

This choice leaves the slack needed in both discretizations.  The first
derivatives of $f$ are affine and hence Lipschitz on a fixed
neighbourhood of $[0,1]^3$.  Replacing each variable by a point at
distance $O(\delta)$ changes $f(u)-f(v)$ by $O_f(\delta)$ and each
$D_j$ by $O_f(\delta)$.  After decreasing $\delta_0$, we have
\[
 C_f\delta\leq\frac{\rho}{2}.
\]
Thus, $|D_j(u,v)|>\rho$ implies that the corresponding derivative at
the representatives is at least $\rho/2$, and the coincidence window
enlarges only to $C_f'\delta$.

\subsubsection{The transverse contribution}
\label{subsubsec:transverse-contribution}

We first record two elementary facts about the line families that arise
below.  The first one explains precisely where the non-degeneracy of
$f$ is used.  Put $v_1=x$, $v_2=y$, and $v_3=z$.  For each
$q\in\{1,2,3\}$, write
\begin{equation}\label{eq:transverse-distinguished-decomposition}
 f(v_1,v_2,v_3)
 =
 A_qv_q^2+B_qv_q
 +m_q(v_p,v_s)v_q+k_q(v_p,v_s),
 \qquad
 \{p,q,s\}=\{1,2,3\},
\end{equation}
where $m_q$ is linear and $k_q$ has degree at most two.  Define
\begin{equation}\label{eq:transverse-Jq}
 J_q(v_p,v_s)
 :=
 \det
 \frac{\partial(m_q,k_q)}{\partial(v_p,v_s)}.
\end{equation}
Thus, $J_q$ is an affine-linear polynomial.

\begin{lemma}
\label{lem:transverse-nondegenerate-J}
If $f$ is non-degenerate as defined in the introduction, then
$J_q\not\equiv0$ for at least one $q\in\{1,2,3\}$.
\end{lemma}

\begin{proof}
After a vertical shifting, we may assume $f(0,0,0)=0$ and write
\[
 f(x,y,z)
 =
 axy+bxz+cyz+dx^2+ey^2+gz^2+hx+iy+jz.
\]
A direct calculation gives
\begin{align*}
 J_1(y,z)
 &=
 (ac-2be)y+(2ag-bc)z+(aj-bi),\\
 J_2(x,z)
 &=
 (ab-2cd)x+(2ag-bc)z+(aj-ch),\\
 J_3(x,y)
 &=
 (ab-2cd)x+(2be-ac)y+(bi-ch),
\end{align*}
up to an inessential simultaneous change of sign in any one of these
expressions.  Suppose, for contradiction, that all three polynomials
vanish identically:
\begin{equation}\label{eq:transverse-coefficient-relations}
 ab=2cd,\qquad
 ac=2be,\qquad
 bc=2ag,
 \qquad
 aj=bi=ch.
\end{equation}
We divide the proof into the following cases.
\begin{enumerate}
\item
If $a=b=c=0$, then
\[
 f(x,y,z)
 =
 (dx^2+hx)+(ey^2+iy)+(gz^2+jz),
\]
which is a forbidden form with outer polynomial equal to the identity.

\item
Suppose that exactly one of $a,b,c$ is nonzero.  After permuting the
variables, we may assume that $a\neq0$ and $b=c=0$.  Then
$bc=2ag$ and $aj=bi$ give $g=j=0$, so $f$ does not depend on $z$,
a contradiction.

\item
If exactly two of $a,b,c$ are nonzero, then, again after a permutation,
$a,b\neq0$ and $c=0$; the identity $ab=2cd$ is impossible.

\item
It remains to consider $abc\neq0$.  In that case,
\[
 d=\frac{ab}{2c},
 \qquad
 e=\frac{ac}{2b},
 \qquad
 g=\frac{bc}{2a}.
\]
In particular, $d\neq0$.  With
\[
 L(x,y,z)
 :=
 x+\frac{a}{2d}y+\frac{b}{2d}z,
\]
the quadratic part of $f$ is $dL^2$.  Moreover,
\[
 \frac{ha}{2d}
 =
 \frac{hc}{b}
 =
 \frac{aj}{b}
 =
 i,
 \qquad
 \frac{hb}{2d}
 =
 \frac{hc}{a}
 =
 j,
\]
where we used $2d=ab/c$ and the last relations in
\eqref{eq:transverse-coefficient-relations}.  Thus,
\[
 hx+iy+jz=hL.
\]
Consequently,
\[
 f(x,y,z)
 =
 dL(x,y,z)^2+hL(x,y,z),
\]
which is again excluded by the definition of non-degenerate quadratic polynomials. The assumption that
all three $J_q$ vanish therefore contradicts non-degeneracy.
\end{enumerate}
\end{proof}

We use the affine-line metric $d_{\mathcal A}$ fixed in
\eqref{eq:inc-line-metric}.  The next lemma replaces the separation
argument based on an auxiliary bridge point.  Its proof integrates only
along the segment joining the two original parameter points; convexity
therefore guarantees that every point used in the proof remains in the
given parameter region.

\begin{lemma}
\label{lem:transverse-line-separation}
Let $m:\mathbb R^2\to\mathbb R$ be a nonconstant linear function, let
$k:\mathbb R^2\to\mathbb R$ have degree at most two, and set
\[
 \ell_u
 :=
 \{(Y,X)\in\mathbb R^2:X=m(u)Y+k(u)\},
 \qquad
 J(u)
 :=
 \det D(m,k)(u).
\]
Suppose that $m$ and $k$ range in a fixed bounded set on a convex set
$\Omega\subset[0,1]^2$.  If $0<\rho\leq1$ and
\[
 |J(u)|\geq\rho
 \qquad
 \text{for all }u\in\Omega,
\]
then
\begin{equation}\label{eq:transverse-line-separation}
 d_{\mathcal A}(\ell_u,\ell_{\widetilde u})
 \geq
 c_{m,k}\rho\,|u-\widetilde u|
 \qquad
 \text{for all }u,\widetilde u\in\Omega.
\end{equation}
In particular, a $\delta$-separated subset of $\Omega$ determines a
$c_{m,k}\rho\delta$-separated family of lines.
\end{lemma}

\begin{proof}
Let
\[
 \lambda:=|\nabla m|>0,
\]
and introduce orthonormal coordinates
\[
 u=se_1+te_2,
\]
where
\[
 e_1=\frac{\nabla m}{\lambda},
 \qquad
 e_2
 =
 \frac{(-\partial_2m,\partial_1m)}{\lambda}.
\]
We have
\begin{equation}\label{eq:m_and_lambda}
 m(se_1+te_2)
 =
 m(0)+\lambda s,
\end{equation}
where the additive constant plays no role after taking differences.
Additionally,
\begin{equation}\label{eq:transverse-kt-J}
 \partial_tk(se_1+te_2)
 =
 \frac{J(se_1+te_2)}{\lambda}.
\end{equation}

Fix $u,\widetilde u\in\Omega$, write their coordinates as
$(s,t)$ and $(\widetilde s,\widetilde t)$, and put
\[
 r:=|u-\widetilde u|,
 \qquad
 \Delta s:=\widetilde s-s,
 \qquad
 \Delta t:=\widetilde t-t.
\]
The case $r=0$ is immediate.  Since $J$ is continuous, never vanishes
on the convex set $\Omega$, and $\Omega$ is connected, $J$ has one
sign throughout $\Omega$.

Let
\[
 C_k:=\sup_{[0,1]^2}|\nabla k|,
\]
and choose a sufficiently small constant $c_{m,k}>0$, depending only
on $m$ and $k$.  We divide the proof into two cases.
\begin{enumerate}
\item
If
\[
 |\Delta s|
 \geq
 c_{m,k}\rho r,
\]
then, by \eqref{eq:m_and_lambda},
\[
 |m(\widetilde u)-m(u)|
 =
 \lambda|\Delta s|
 \geq
 c_{m,k}\lambda\rho r.
\]

\item
Suppose instead that
\[
 |\Delta s|
 <
 c_{m,k}\rho r.
\]
Taking $c_{m,k}\leq1/2$, we have $|\Delta t|\geq r/2$.  The segment
\[
 u_\tau
 :=
 (1-\tau)u+\tau\widetilde u,
 \qquad
 0\leq\tau\leq1,
\]
is contained in $\Omega$.  Hence, by the fundamental theorem of
calculus,
\[
 k(\widetilde u)-k(u)
 =
 \Delta s\int_0^1\partial_sk(u_\tau)\,d\tau
 +
 \Delta t\int_0^1\partial_tk(u_\tau)\,d\tau.
\]
The second integral has no cancellation because $J$ has constant sign.
Together with \eqref{eq:transverse-kt-J} and the assumption
$|J(u)|\geq\rho$, we obtain
\[
\begin{aligned}
 |k(\widetilde u)-k(u)|
 &\geq
 \frac{\rho}{\lambda}|\Delta t|
 -
 C_k|\Delta s|\\
 &\geq
 \left(
 \frac{1}{2\lambda}-C_kc_{m,k}
 \right)\rho r\\
 &\gtrsim_{m,k}
 \rho r,
\end{aligned}
\]
after decreasing $c_{m,k}$ if necessary.
\end{enumerate}
Hence, in both cases,
\begin{equation}\label{eq:maxmk_lb}
 \max\bigl\{
 |m(\widetilde u)-m(u)|,
 |k(\widetilde u)-k(u)|
 \bigr\}
 \gtrsim_{m,k}
 \rho r.
\end{equation}

It remains only to compare slope-intercept coordinates with
$d_{\mathcal A}$.  On the bounded range in question, the line
$X=mY+k$ has the canonical representative
\[
 \omega
 =
 \frac{(-m,1)}{\sqrt{1+m^2}},
 \qquad
 b
 =
 \frac{k}{\sqrt{1+m^2}}.
\]
Here the second coordinate of $\omega$ is bounded below by a positive
constant.  The inverse formulas
\[
 m=-\frac{\omega_1}{\omega_2},
 \qquad
 k=\frac{b}{\omega_2},
\]
show that the map $(m,k)\mapsto\ell_{m,k}$ and its inverse are
Lipschitz on this compact set.  The sign-reversed canonical
representatives remain a fixed positive distance apart; since the
parameter range is bounded, taking the quotient by this sign does not
alter the lower Lipschitz comparison.  Thus, combining this observation
with \eqref{eq:maxmk_lb}, we obtain
\[
 d_{\mathcal A}(\ell_u,\ell_{\widetilde u})
 \gtrsim_{m,k}
 |m(u)-m(\widetilde u)|
 +
 |k(u)-k(\widetilde u)|,
\]
which completes the proof of
\eqref{eq:transverse-line-separation}.
\end{proof}

We now give the discrete statement needed both for AD-regular measures
and for the dyadic reduction of a general Frostman measure.  Allowing
six different sets, and assuming only upper cardinality bounds, is
important in the latter application.

\begin{proposition}
\label{prop:six-set-transverse}
Let $f\in\mathbb R[x,y,z]$ be a fixed non-degenerate polynomial of
degree at most two.  Let $0<\alpha<1$, let
\[
 C_0,C_1,c_1^{-1}\geq1,
\]
and put
\[
 \sigma:=\min\{2\alpha,2-2\alpha\}.
\]
Choose
\begin{equation}\label{eq:transverse-kappa-choice}
 0<\kappa<\frac{\alpha^2\sigma}{64}.
\end{equation}
Then there are
\[
 \eta=\eta(\alpha,\kappa)>0,
\]
\[
 \delta_0
 =
 \delta_0(f,\alpha,\kappa,C_0,C_1,c_1)>0,
\]
and
\[
 C
 =
 C(f,\alpha,\kappa,C_0,C_1,c_1)<\infty
\]
such that the following holds for $0<\delta\leq\delta_0$.

For $i=1,2,3$, let
\[
 M_i,M_i'\subset[0,1]
\]
be $\delta$-separated sets satisfying
\begin{equation}\label{eq:six-set-NC}
 |M_i|,|M_i'|
 \leq
 C_0\delta^{-\alpha},
 \qquad
 |M_i\cap I|,|M_i'\cap I|
 \leq
 C_0|I|^\alpha\delta^{-\alpha}
\end{equation}
for every interval $I$ with $|I|\geq\delta$.  For every
$j\in\{1,2,3\}$, one has
\begin{align}
 &\#\Bigl\{
 (\vec v,\vec v')
 \in
 (M_1\times M_2\times M_3)
 \times
 (M_1'\times M_2'\times M_3'):
 \notag\\[-2mm]
 &\hspace{24mm}
 |f(\vec v)-f(\vec v')|<
 C_1\delta,
 \quad
 |\partial_jf(\vec v')|
 \geq
 c_1\delta^\kappa
 \Bigr\}
 \leq
 C\delta^{-5\alpha+\eta}.
 \label{eq:six-set-transverse-count}
\end{align}
The same conclusion holds when the derivative condition is imposed on
$\vec v$ rather than on $\vec v'$.
\end{proposition}

\begin{proof}
We prove \eqref{eq:six-set-transverse-count}; the last assertion follows
by interchanging the primed and unprimed triples.  By
Lemma~\ref{lem:transverse-nondegenerate-J}, we may choose a
distinguished unprimed coordinate $q\in\{1,2,3\}$ for which
$J_q\not\equiv0$.  Fix an arbitrary $j\in\{1,2,3\}$.  Let $p,s$
denote the remaining two unprimed coordinate indices:
\[
 \{p,s\}
 =
 \{1,2,3\}\setminus\{q\}.
\]

Fix the two primed coordinates $(v_i')_{i\neq j}$.  As a function of
the remaining variable, write
\[
 h(t)
 :=
 \begin{cases}
 f(t,v_2',v_3'),
 &\text{if }j=1,\\
 f(v_1',t,v_3'),
 &\text{if }j=2,\\
 f(v_1',v_2',t),
 &\text{if }j=3.
 \end{cases}
\]
Then
\[
 \partial_jf(\vec v')
 =
 h'(t)
 =
 2D_jt+L,
\]
where $D_j$ is the coefficient of $v_j^2$ in $f$, and $L$ depends on
$j$ and on the two fixed primed coordinates.  We divide the proof into
two cases.
\begin{enumerate}
\item
If $D_j\neq0$, the set
\begin{equation}\label{eq:set_with_lb}
 \bigl\{
 t\in[0,1]:
 |h'(t)|\geq c_1\delta^\kappa
 \bigr\}
\end{equation}
is the union of at most two intervals, and $h'(t)$ has constant sign
and satisfies
\[
 |h'(t)|
 \geq
 c_1\delta^\kappa
\]
on each interval.

\item
If $D_j=0$, then $h'(t)$ is constant.  The set in
\eqref{eq:set_with_lb} is therefore either all of $[0,1]$ or empty.
If the set is empty, then \eqref{eq:six-set-transverse-count} holds
trivially.
\end{enumerate}
Thus, we may assume that the set in \eqref{eq:set_with_lb} is nonempty.
Moreover, to prove \eqref{eq:six-set-transverse-count}, it is enough to
work on one such interval $A$ and replace $M_j'$ by $M_j'\cap A$.

Using the decomposition
\eqref{eq:transverse-distinguished-decomposition}, define
\begin{equation}\label{eq:transverse-point-map}
 \Phi(t,w)
 :=
 \bigl(
 h(t)-A_qw^2-B_qw,w
 \bigr),
 \qquad
 (t,w)\in(M_j'\cap A)\times M_q.
\end{equation}
This map is bi-Lipschitz with the precise scale dependence needed in
Theorem~\ref{our_incidence_theorem}.  Indeed, the upper Lipschitz bound
is uniform because $f$ is fixed and all variables lie in $[0,1]$.
For two points $(t,w)$ and $(\widetilde t,\widetilde w)$ in its
domain, the mean value theorem on the interval $A$ gives
\[
 c_1\delta^\kappa|t-\widetilde t|
 \leq
 |h(t)-h(\widetilde t)|.
\]
On the other hand, for some constant $C_f'>0$ depending only on $f$,
\[
 |h(t)-h(\widetilde t)|
 \leq
 |\Phi_1(t,w)-\Phi_1(\widetilde t,\widetilde w)|
 +
 C_f'|w-\widetilde w|.
\]
Since
\[
 w=\Phi_2(t,w),
\]
these two inequalities imply that there exists a constant
$C_{f,c_1}$, depending only on $f$ and $c_1$, such that
\begin{equation}\label{eq:transverse-Phi-bilip}
 C_{f,c_1}^{-1}
 \delta^\kappa
 |(t,w)-(\widetilde t,\widetilde w)|
 \leq
 |\Phi(t,w)-\Phi(\widetilde t,\widetilde w)|
 \leq
 C_{f,c_1}
 \delta^{-\kappa}
 |(t,w)-(\widetilde t,\widetilde w)|.
\end{equation}
The image of $\Phi$ lies in a fixed bounded set.

For $(u,v)\in M_p\times M_s$, let
\begin{equation}\label{eq:transverse-lines}
 \ell_{u,v}:
 \quad
 X=m_q(u,v)Y+k_q(u,v).
\end{equation}
By construction, if
\[
 P_{t,w}:=\Phi(t,w),
\]
then
\begin{equation}\label{eq:coincidence-is-incidence}
 \bigl|f(\vec v)-f(\vec v')\bigr|
 =
 \bigl|
 P_{t,w}^{(1)}
 -
 m_q(u,v)P_{t,w}^{(2)}
 -
 k_q(u,v)
 \bigr|.
\end{equation}
All slopes are uniformly bounded.  Consequently, the condition that
the left-hand side of \eqref{eq:coincidence-is-incidence} be at most
$C_1\delta$ places $P_{t,w}$ in the $C_f\delta$-neighbourhood of
$\ell_{u,v}$.  This is the required enlarged coincidence window; no
shrinking from $C_f\delta$ to $\delta$ is used.

We next separate the parameters for which the affine Jacobian is small.
Let
\[
 \Omega_{\mathrm{bad}}
 :=
 \bigl\{
 (u,v)\in M_p\times M_s:
 |J_q(u,v)|<\delta^\gamma
 \bigr\},
\]
where $0<\gamma<1$ will be chosen below.  If $J_q$ is a nonzero
constant, then
\[
 \Omega_{\mathrm{bad}}=\varnothing
\]
for all sufficiently small $\delta$.  Otherwise, at least one of its
two linear coefficients is nonzero.  After fixing the other coordinate,
the corresponding section of
\[
 \{|J_q|<\delta^\gamma\}
\]
is an interval of length $O_f(\delta^\gamma)$.  Since
$\delta^\gamma\geq\delta$, \eqref{eq:six-set-NC} gives
\begin{equation}\label{eq:bad-J-cardinality}
 |\Omega_{\mathrm{bad}}|
 \lesssim_{f,C_0}
 \delta^{-2\alpha+\alpha\gamma}.
\end{equation}

For a fixed $(u,v)\in\Omega_{\mathrm{bad}}$ and a fixed $w\in M_q$,
the condition
\[
 |h(t)-f(\vec v)|<C_1\delta
\]
limits $h(t)$ to an interval of length $O_{f,C_1}(\delta)$.  On the
interval $A$, the function $h$ is monotone and
\[
 |h'|\geq c_1\delta^\kappa.
\]
Hence, the admissible values of $t$ lie in an interval of length
\[
 O_{f,c_1,C_1}(\delta^{1-\kappa}),
\]
which contains at most
\[
 O_{f,c_1,C_1}(\delta^{-\kappa})
\]
points of the $\delta$-separated set $M_j'\cap A$.  Summing over
$w\in M_q$ and using \eqref{eq:bad-J-cardinality}, we obtain
\begin{equation}\label{eq:bad-J-count}
\begin{aligned}
 &\#\Bigl\{
 (t,v_p,v_s,v_q)
 \in
 (M_j'\cap A)\times M_p\times M_s\times M_q:
 \\
 &\hspace{27mm}
 (v_p,v_s)\in\Omega_{\mathrm{bad}},
 \quad
 |h(t)-f(\vec v)|<C_1\delta
 \Bigr\}
 \\
 &\hspace{12mm}
 \lesssim_{f,c_1,C_1,C_0}
 \delta^{-3\alpha+\alpha\gamma-\kappa},
\end{aligned}
\end{equation}
where
\[
 (v_p,v_s,v_q)=(u,v,w),
\]
and $\vec v=(v_1,v_2,v_3)$ equals $(v_p,v_s,v_q)$ up to a coordinate
permutation.

The complement of the strip is the union of the two convex regions
\begin{align*}
 \Omega_+
 &:=
 \bigl\{
 (u,v)\in[0,1]^2:
 J_q(u,v)\geq\delta^\gamma
 \bigr\},\\
 \Omega_-
 &:=
 \bigl\{
 (u,v)\in[0,1]^2:
 J_q(u,v)\leq-\delta^\gamma
 \bigr\}.
\end{align*}
Either region is allowed to be empty; when $J_q$ is constant, one of
them is all of $[0,1]^2$ for small $\delta$.  The parameter set
$M_p\times M_s$ is $\delta$-separated.
Lemma~\ref{lem:transverse-line-separation}, applied to each of
$\Omega_+$ and $\Omega_-$ with $\rho=\delta^\gamma$, shows that the
corresponding lines form a
\[
 c_f\delta^{1+\gamma}
\]
separated family in $d_{\mathcal A}$.

For completeness, such a family can be partitioned into
$O_f(\delta^{-2\gamma})$ many $c_f'\delta$-separated subfamilies.
Indeed, join two lines when their $d_{\mathcal A}$-distance is less
than $c_f'\delta$.  The lines lie in a fixed compact subset of the
two-dimensional affine-line space.  A packing estimate in a line-space
ball of radius $O_f(\delta)$ shows that every vertex has degree
\[
 O_f\left(
 \left(
 \frac{\delta}{\delta^{1+\gamma}}
 \right)^2
 \right).
\]
Greedy colouring gives the stated partition.  The fixed-constant
robustness of Theorem~\ref{our_incidence_theorem} allows both
$c_f'\delta$-separation and $C_f\delta$-neighbourhoods; equivalently,
one may replace $\delta$ by a fixed multiple and apply a bounded
additional colouring.

The two variables $(t,w)$ in \eqref{eq:transverse-point-map} satisfy
\eqref{eq:six-set-NC} and have cardinality at most
$C_0\delta^{-\alpha}$.  The point map $\Phi$ satisfies
\eqref{eq:transverse-Phi-bilip}, while each line subfamily has
cardinality at most $C_0^2\delta^{-2\alpha}$.  The rectangular,
upper-cardinality form of Theorem~\ref{our_incidence_theorem} therefore
applies to every colour class, after enlarging its structural constant
by a fixed factor depending only on $f,C_0,C_1,c_1$.  If
$\varepsilon_{\mathrm{inc}}$ denotes its gain in
Theorem~\ref{our_incidence_theorem} for each subfamily, then the total
contribution of $\Omega_+\cup\Omega_-$, for the fixed primed
coordinates and the fixed derivative branch, is
\begin{equation}\label{eq:good-J-count}
\begin{aligned}
 &\#\Bigl\{
 (t,v_p,v_s,v_q)
 \in
 (M_j'\cap A)\times M_p\times M_s\times M_q:
 \\
 &\hspace{25mm}
 (v_p,v_s)\notin\Omega_{\mathrm{bad}},
 \quad
 |h(t)-f(\vec v)|<C_1\delta
 \Bigr\}
 \\
 &\hspace{12mm}
 \lesssim_{f,c_1,C_1,C_0}
 \delta^{-3\alpha+\varepsilon_{\mathrm{inc}}-2\gamma}.
\end{aligned}
\end{equation}
The extra $2\gamma$ in the exponent comes from the
$O_f(\delta^{-2\gamma})$ subfamilies.

We now make the parameter hierarchy explicit.  By
\eqref{eq:inc-epsilon-range}, we may choose any positive
$\varepsilon_{\mathrm{inc}}<G/4$.  Let
\begin{equation}\label{eq:transverse-parameter-hierarchy}
 G
 :=
 \alpha\left(
 \frac{\sigma}{2}-\kappa
 \right),
 \qquad
 \varepsilon_{\mathrm{inc}}
 :=
 \frac{\alpha\sigma}{16},
 \qquad
 \gamma
 :=
 \frac{\alpha\sigma}{48}.
\end{equation}
The choice \eqref{eq:transverse-kappa-choice} gives
\[
 \frac{2\kappa}{\alpha}
 <
 \varepsilon_{\mathrm{inc}}
 <
 \frac G4,
 \qquad
 \frac{\kappa}{\alpha}
 <
 \gamma
 <
 \frac{\varepsilon_{\mathrm{inc}}}{2}.
\]
Thus, Theorem~\ref{our_incidence_theorem} is available with the chosen
gain, and
\begin{equation}\label{eq:transverse-positive-gain}
 \eta_0
 :=
 \min\{
 \varepsilon_{\mathrm{inc}}-2\gamma,
 \alpha\gamma-\kappa
 \}
 >0.
\end{equation}
For example, the choices
\eqref{eq:transverse-kappa-choice} and
\eqref{eq:transverse-parameter-hierarchy} imply
\[
 \eta_0>\frac{\alpha^2\sigma}{192},
\]
so one may take
\[
 \eta=\frac{\alpha^2\sigma}{384}.
\]

Combining \eqref{eq:bad-J-count} and \eqref{eq:good-J-count}, and
summing over the at most two derivative branches, gives
\[
\begin{aligned}
 &\#\Bigl\{
 (v_j',v_q,v_p,v_s)
 \in
 (M_j'\cap A)\times M_q\times M_p\times M_s:
 \\
 &\hspace{25mm}
 |h(v_j')-f(\vec v)|<C_1\delta,
 \quad
 |h'(v_j')|\geq c_1\delta^\kappa
 \Bigr\}
 \\
 &\hspace{12mm}
 \lesssim_{f,c_1,C_1,C_0}
 \delta^{-3\alpha+\eta}
\end{aligned}
\]
for each fixed pair $(v_i')_{i\neq j}$.  Finally, there are at most
\[
 C_0^2\delta^{-2\alpha}
\]
choices for that pair.  This proves
\eqref{eq:six-set-transverse-count}.
\end{proof}

We next explain how Proposition~\ref{prop:six-set-transverse} enters
the measure estimates.  This also records the two discretization
details that will be used below: the coincidence window is enlarged by
a fixed constant, and the derivative threshold is retained with a
fixed amount of slack.

\begin{lemma}\label{lem:AD-transverse-bound}
Let $\mu$ be an $\alpha$-AD-regular probability measure supported on
$[0,1]$, and choose $\kappa$ as in
\eqref{eq:transverse-kappa-choice}.  There exists
$\eta=\eta(\alpha,\kappa)>0$ such that
\[
 I_i
 \lesssim
 \delta^{\alpha+\eta},
 \qquad
 1\leq i\leq6.
\]
\end{lemma}

\begin{proof}
Recall the definition of $I_j$:
\[
 I_j
 :=
 \iint
 \mathbf 1_{\{|f(u)-f(v)|\leq2\delta\}}
 \mathbf 1_{\{|D_j(u,v)|>\rho\}}
 \,d\mu_3(u)\,d\mu_3(v),
 \qquad
 1\leq j\leq6.
\]

We give the detailed proof for $I_1$ first and then explain how the same
method applies to $I_2,\ldots,I_6$.  Let $M$ be a maximal
$\delta$-separated subset of $\operatorname{spt}\mu$.  AD-regularity
and the usual disjoint-ball argument give
\[
 |M|
 \lesssim
 \delta^{-\alpha},
 \qquad
 |M\cap J|
 \lesssim
 |J|^\alpha\delta^{-\alpha}
 \quad
 (|J|\geq\delta).
\]
Associate to each point of the support a nearest point of $M$.  Every
resulting cell has $\mu$-mass $O(\delta^\alpha)$.  Since $f$ and its
first derivatives have bounded derivatives on a fixed neighbourhood of
$[0,1]^3$, replacing the six variables by their cell representatives
changes $f(\vec v)-f(\vec v')$ by at most $C_f\delta$ and changes
$\partial_xf(\vec v')$ by at most $C_f\delta$.  For sufficiently small
$\delta$, the implication
\[
 |\partial_xf(\vec v')|
 >
 \delta^\kappa
 \quad\Longrightarrow\quad
 |\partial_xf(\widetilde{\vec v'})|
 \geq
 \frac12\delta^\kappa
\]
therefore holds for the representative $\widetilde{\vec v'}$.
Proposition~\ref{prop:six-set-transverse}, with $C_1=C_f$ and
$c_1=1/2$, bounds the number of representative six-tuples by
\[
 O(\delta^{-5\alpha+\eta}).
\]
Multiplying by the $O(\delta^{6\alpha})$ upper bound for the product of
the six cell masses gives
\[
 I_1
 \lesssim
 \delta^{\alpha+\eta}.
\]
The same argument applies to $I_2$ and $I_3$.  For $I_4,I_5,I_6$,
interchange the primed and unprimed triples and use the last assertion
of Proposition~\ref{prop:six-set-transverse}.
\end{proof}

The AD-regular assumption has been used above only in the unweighted
discretization of the transverse contribution.  The small-gradient
analysis below requires only the upper Frostman bound and will be used
unchanged in the proof for general Frostman measures in
Subsection~\ref{sec:Frostman}.

\subsubsection{Bounding \texorpdfstring{$I_0$}{}}

Recall that
\[
 I_0
 :=
 \iint
 \mathbf 1_{\{|f(u)-f(v)|\leq2\delta\}}
 \mathbf 1_{\mathcal G_{2\rho}}(u)
 \mathbf 1_{\mathcal G_{2\rho}}(v)
 \,d\mu_3(u)\,d\mu_3(v),
\]
where
\[
 \mathcal G_s
 :=
 \bigl\{
 u\in\operatorname{spt}(\mu)^3:
 |\nabla f(u)|\leq s
 \bigr\}.
\]

\begin{lemma}\label{lem:I0_bound}
Let $I_0$ be defined by \eqref{eq:def_I0}.  There is
\[
 \delta_0
 =
 \delta_0(\alpha,\kappa,f,C_\mu)>0
\]
such that
\[
 I_0
 \lesssim
 \delta^{\alpha+3\alpha\kappa/2}
 \qquad
 (0<\delta\leq\delta_0).
\]
Here $\mu$ need only be an $\alpha$-Frostman probability measure
supported on $[0,1]$; Ahlfors--David regularity is not required.
\end{lemma}

To prove this lemma, we first observe some geometric properties of the
small-gradient set.

Write
\[
 u=(x,y,z)\in\mathbb R^3.
\]
Since $f$ is quadratic, there exist a symmetric matrix
$H\in\mathbb R^{3\times3}$, a vector $\vec b\in\mathbb R^3$, and a
scalar $c_0\in\mathbb R$ such that
\begin{equation}\label{eq:quad_form}
 f(u)
 =
 \frac12u^THu+\vec b\cdot u+c_0,
 \qquad
 \nabla f(u)
 =
 Hu+\vec b.
\end{equation}
Let
\[
 K
 :=
 \bigl\{
 u\in\mathbb R^3:
 \nabla f(u)=0
 \bigr\}
 =
 \bigl\{
 u\in\mathbb R^3:
 Hu+\vec b=0
 \bigr\}.
\]
If $K=\emptyset$, then $|\nabla f|$ is bounded below on $[0,1]^3$ by
a positive constant depending only on $f$.  In this case,
\[
 \mathcal G_{2\delta^\kappa}=\emptyset
\]
for all sufficiently small $\delta$, and hence $I_0=0$.  Thus, we may
assume that $K\neq\emptyset$.

\begin{lemma}\label{lem:K_dimension}
If $K\neq\emptyset$, then
\[
 \operatorname{rank}H\in\{2,3\}.
\]
Consequently, $K$ is either an affine line or a single point.
\end{lemma}

\begin{proof}
Choose $u_0\in K$.  Translating to $u_0$ gives
\[
 f(u)
 =
 f(u_0)
 +
 \frac12(u-u_0)^TH(u-u_0).
\]
If $H=0$, then $f$ is constant, contrary to the assumption that it
depends on each variable.  If $\operatorname{rank}H=1$, then
\[
 H=\lambda vv^T
\]
for some $\lambda\neq0$ and
$v\in\mathbb R^3\setminus\{0\}$.  Hence
\[
 f(x,y,z)
 =
 f(u_0)
 +
 \frac{\lambda}{2}
 \bigl(
 v_1(x-x_0)+v_2(y-y_0)+v_3(z-z_0)
 \bigr)^2,
\]
which has the forbidden form $G(I(x)+J(y)+K(z))$.  Thus,
\[
 \operatorname{rank}H\geq2.
\]
Finally, by construction,
\[
 K=u_0+\ker H.
\]
Hence, $K$ is either an affine line or a single point.
\end{proof}

\begin{lemma}\label{lem:Gr_tube}
If $K\neq\emptyset$, there exists a constant $C\geq1$, depending only
on $f$, such that, for all sufficiently small $r>0$,
\[
 \mathcal G_r
 \subset
 N_{Cr}(K)\cap[0,1]^3,
\]
where $N_\varrho(K)$ denotes the Euclidean $\varrho$-neighbourhood of
$K$.
\end{lemma}

\begin{proof}
Fix $u_0\in K$.  By \eqref{eq:quad_form},
\[
 \nabla f(u)
 =
 H(u-u_0)
\]
for all $u$.  Let
\[
 V:=(\ker H)^\perp.
\]
The restriction $H|_V$ is invertible, since $V$ is the range of $H$.
Let $\sigma_{\min}$ denote the smallest singular value of $H|_V$.
Then
\[
 \sigma_{\min}>0.
\]

For any $u\in\mathbb R^3$, write
\[
 u-u_0=v+w,
 \qquad
 v\in V,
 \quad
 w\in\ker H.
\]
Then
\[
 H(u-u_0)=Hv,
\]
and hence
\[
 |\nabla f(u)|
 =
 |Hv|
 \geq
 \sigma_{\min}|v|.
\]
If $|\nabla f(u)|\leq r$, then
\[
 |v|
 \leq
 \sigma_{\min}^{-1}r.
\]
Since $u_0+w\in K$ and $v\perp\ker H$, the distance from $u$ to $K$
equals $|v|$.  Thus,
\[
 u\in N_{\sigma_{\min}^{-1}r}(K).
\]
This proves the lemma with
\[
 C=\sigma_{\min}^{-1}.
\]
\end{proof}

The next lemma bounds the measure of tubes under $\mu^{\otimes3}$.

\begin{lemma}\label{lem:tube_mass_line}
Let $\mu$ be an $\alpha$-Frostman probability measure supported on
$[0,1]$.  Let $\Lambda\subset\mathbb R^3$ be an affine line.  Then,
for all $0<r\leq1$,
\begin{equation}\label{eq:tube-mass-line}
 \mu^{\otimes3}\bigl(N_r(\Lambda)\bigr)
 \lesssim
 r^{2\alpha}.
\end{equation}
If $p\in\mathbb R^3$, then, for all $0<r\leq1$,
\begin{equation}\label{eq:point-mass-ball}
 \mu^{\otimes3}\bigl(B(p,r)\bigr)
 \lesssim
 r^{3\alpha}.
\end{equation}
\end{lemma}

\begin{proof}
The point estimate \eqref{eq:point-mass-ball} follows directly from the
definition of a Frostman measure.

We now prove the line estimate \eqref{eq:tube-mass-line}.  Let
$\Lambda$ be an affine line.  After permuting coordinates, $\Lambda$
has a parametrization of the form
\[
 \Lambda
 =
 \bigl\{
 (t,\phi(t),\psi(t)):
 t\in\mathbb R
 \bigr\},
\]
where $\phi$ and $\psi$ are affine-linear functions.  Since $[0,1]^3$
is compact, there is a constant $C\geq1$, depending only on $\Lambda$,
such that
\[
 N_r(\Lambda)\cap[0,1]^3
 \subset
 \bigl\{
 (t,y,z)\in[0,1]^3:
 |y-\phi(t)|\leq Cr,\ 
 |z-\psi(t)|\leq Cr
 \bigr\}.
\]
Therefore,
\begin{equation}\label{eq:tube_mass_fubini}
 \mu^{\otimes3}\bigl(N_r(\Lambda)\bigr)
 \leq
 \int
 \mu\bigl(B(\phi(t),Cr)\bigr)
 \mu\bigl(B(\psi(t),Cr)\bigr)
 \,d\mu(t).
\end{equation}
Since $\mu$ is an $\alpha$-Frostman measure,
\[
 \mu(B(\cdot,Cr))
 \lesssim
 r^\alpha
\]
uniformly.  Hence,
\[
 \mu^{\otimes3}\bigl(N_r(\Lambda)\bigr)
 \lesssim
 r^{2\alpha}.
\]
\end{proof}

We next control the mass of a $\delta$-thick level set of a rank-two
quadratic form under $\mu\times\mu$.

\begin{lemma}\label{lem:rank2_sublevel}
Let \(0<\alpha\leq1\), and let \(\mu\) be an
\(\alpha\)-Frostman probability measure supported on \([0,1]\):
\begin{equation}
 \mu(B(x,r))\leq C_\mu r^\alpha
 \qquad(x\in\mathbb{R},\ r>0).
 \label{eq:frostman-uniform}
\end{equation}
Fix \(M\geq1\) and \(\kappa>0\).  Consider any real quadratic
polynomial
\begin{equation}
 Q(u,v)=au^2+buv+cv^2+du+ev+g
 \label{eq:quadratic-uniform}
\end{equation}
satisfying
\begin{equation}
 \max\{|a|,|b|,|c|,|d|,|e|,|g|\}\leq M,
 \qquad
 |4ac-b^2|\geq\kappa.
 \label{eq:coefficient-hessian-bounds}
\end{equation}
Then there is a constant
\[
 C=C(\alpha,\mu,M,\kappa)<\infty
\]
such that, for every \(t\in\mathbb{R}\) and sufficientlly small $\delta>0$,
\begin{equation}
 (\mu\times\mu)
 \bigl\{(u,v)\in[0,1]^2:
       |Q(u,v)-t|\leq\delta\bigr\}
 \leq
 C\delta^\alpha\log\left(\frac1\delta\right).
 \label{eq:uniform-rank2-conclusion}
\end{equation}
In particular, the constant is uniform over every family with a common
coefficient bound \(M\) and a common positive lower bound \(\kappa\)
for the absolute Hessian determinant.
\end{lemma}

\begin{proof}
We first record two elementary one-dimensional estimates.

\medskip
\noindent
\textbf{Claim 1: sublevel intervals of a one-variable quadratic.}
Consider a quadratic polynomial
\[p(s)=As^2+Bs+C.\]
For every \(\rho>0\), the set
\[\{s\in\mathbb{R}:|p(s)|\leq\rho\}\]
is a union of at most two intervals, each of length at most
\begin{equation}
 4\sqrt{\rho/|A|}.
 \label{eq:one-variable-quadratic-length}
\end{equation}
Consequently,
\begin{equation}
 \mu\{s\in[0,1]:|p(s)|\leq\rho\}
 \leq
 C_{\alpha,\mu}\left(\frac{\rho}{|A|}\right)^{\alpha/2}.
 \label{eq:one-variable-quadratic-mass}
\end{equation}

Indeed, after completing the square and dividing by \(A\), the relevant
set has the form
\[
 \bigl\{s:|(s-s_0)^2-\gamma|\leq\varepsilon\bigr\},
 \qquad
 \varepsilon:=\rho/|A|.
\]
If \(\gamma\leq-\varepsilon\), the set is empty, apart from a harmless
endpoint when equality holds.  If
\(-\varepsilon<\gamma\leq\varepsilon\), it is contained in
\[
 [s_0-\sqrt{2\varepsilon},s_0+\sqrt{2\varepsilon}].
\]
If \(\gamma>\varepsilon\), it is the union of two intervals, each of
length
\[
 \sqrt{\gamma+\varepsilon}-\sqrt{\gamma-\varepsilon}
 =
 \frac{2\varepsilon}
 {\sqrt{\gamma+\varepsilon}+\sqrt{\gamma-\varepsilon}}
 \leq2\sqrt{\varepsilon}.
\]
This proves \eqref{eq:one-variable-quadratic-length}.
Equation \eqref{eq:one-variable-quadratic-mass} follows from the
Frostman bound, after intersecting these intervals with \([0,1]\).
The estimate remains valid when the displayed upper bound for the
interval length exceeds \(1\), after enlarging the implicit constant,
since \(\mu\) is a probability measure.

\medskip
\noindent
\textbf{Claim 2: a uniform linear-in-\(v\) slicing estimate.}
Fix \(\beta>0\) and \(M_1\geq1\).  Suppose
\[
 \Lambda(u)=bu+e,
 \qquad
 |b|\geq\beta,
 \qquad
 |b|+|e|\leq M_1,
\]
and let \(A:[0,1]\to\mathbb{R}\) be any Borel function.  Then, uniformly
in \(A\),
\begin{equation}
 (\mu\times\mu)
 \bigl\{(u,v)\in[0,1]^2:
       |\Lambda(u)v+A(u)|\leq\rho\bigr\}
 \leq
 C_{\alpha,\mu,\beta,M_1}
 \rho^\alpha\log\left(\frac2\rho\right)
 \label{eq:linear-slicing-estimate}
\end{equation}
for \(0<\rho\leq1/2\).

To prove this, put
\[
 W_{-1}:=\{u\in[0,1]:|\Lambda(u)|\leq\rho\}
\]
and, for \(m\geq0\),
\[
 W_m:=\{u\in[0,1]:
          2^m\rho<|\Lambda(u)|\leq2^{m+1}\rho\}.
\]
The set \(W_{-1}\) is contained in an interval of length
\(2\rho/|b|\), and \(W_m\) is contained in the interval
\(\{|\Lambda|\leq2^{m+1}\rho\}\), whose length is
\(2^{m+2}\rho/|b|\).  Therefore,
\begin{equation}
 \mu(W_{-1})\lesssim_{\alpha,\mu,\beta}\rho^\alpha,
 \qquad
 \mu(W_m)\lesssim_{\alpha,\mu,\beta}
                 (2^m\rho)^\alpha.
 \label{eq:linear-u-shell-mass}
\end{equation}
For fixed \(u\in W_m\), the admissible values of \(v\) in \eqref{eq:linear-slicing-estimate} form an interval
of length at most \(2\rho/|\Lambda(u)|\).  Hence
\begin{equation}
 \mu\{v\in[0,1]:
       |\Lambda(u)v+A(u)|\leq\rho\}
 \lesssim_{\alpha,\mu}2^{-m\alpha}.
 \label{eq:linear-v-fibre-mass}
\end{equation}
The contribution of \(W_{-1}\) is at most
\(C\rho^\alpha\), using the trivial bound \(1\) for the \(v\)-fibre.
Equations \eqref{eq:linear-u-shell-mass} and
\eqref{eq:linear-v-fibre-mass} show that every \(W_m\) contributes at
most \(C\rho^\alpha\).  Since
\[
 |\Lambda(u)|\leq |b|+|e|\leq M_1
 \qquad(u\in[0,1]),
\]
only \(O_{M_1}(\log(2/\rho))\) shells can be nonempty.  This proves
\eqref{eq:linear-slicing-estimate}.

\medskip
\noindent
Now, we are ready to prove Lemma \ref{lem:rank2_sublevel}: Put
\begin{equation}
 \delta_0:=
 \min\left\{\frac14,\frac{\kappa}{8M}\right\}.
 \label{eq:uniform-small-scale-threshold}
\end{equation}
We henceforth assume \(0<\delta\leq\delta_0\).

For \(u\in[0,1]\), write
\begin{equation}
 q_u(v):=Q(u,v)-t
 =cv^2+(bu+e)v+(au^2+du+g-t).
 \label{eq:quadratic-fibre}
\end{equation}
We split according to the size of \(c\) relative to \(\delta\).

\medskip
\noindent
\textbf{Case 1: \(\boldsymbol{|c|\leq\delta}\).}
Remove the small quadratic term in \(v\) and put
\[
 \widetilde q_u(v)
 :=q_u(v)-cv^2
 =(bu+e)v+(au^2+du+g-t).
\]
Since \(0\leq v\leq1\),
\begin{equation}
 |q_u(v)|\leq\delta
 \quad\Longrightarrow\quad
 |\widetilde q_u(v)|\leq2\delta.
 \label{eq:absorb-small-c}
\end{equation}
Moreover, the determinant hypothesis gives
\begin{equation}
 b^2
 \geq |4ac-b^2|-4|a||c|
 \geq\kappa-4M\delta
 \geq\frac{\kappa}{2}.
 \label{eq:mixed-coefficient-lower-bound}
\end{equation}
Thus \(|b|\geq\sqrt{\kappa/2}\), uniformly.  Apply Claim~2 with
\[
 \rho=2\delta,
 \qquad
 \beta=\sqrt{\kappa/2},
 \qquad
 M_1=2M,
 \qquad
 A(u)=au^2+du+g-t.
\]
Since \(2\delta\leq1/2\), \eqref{eq:absorb-small-c} yields
\begin{equation}
 (\mu\times\mu)
 \{(u,v):|q_u(v)|\leq\delta\}
 \lesssim_{\alpha,\mu,M,\kappa}
 \delta^\alpha\log(2/\delta).
 \label{eq:small-c-final}
\end{equation}

\medskip
\noindent
\textbf{Case 2: \(\boldsymbol{|c|\geq\delta}\).}
If the sublevel set is empty, there is nothing to
prove.  Otherwise, for some \((u,v)\in[0,1]^2\),
\(|Q(u,v)-t|\leq\delta\); hence
\begin{equation}
 |t|\leq6M+1.
 \label{eq:uniform-t-bound}
\end{equation}
Define
\[
 B(u):=bu+e,
 \qquad
 A_t(u):=au^2+du+g-t,
\]
and the discriminant
\begin{align}
 \Delta_t(u)
 &:=B(u)^2-4cA_t(u) \notag\\
 &=(b^2-4ac)u^2+(2be-4cd)u
   +(e^2-4c(g-t)).
 \label{eq:uniform-discriminant}
\end{align}
Its leading coefficient has absolute value at least \(\kappa\).
Equations \eqref{eq:coefficient-hessian-bounds} and
\eqref{eq:uniform-t-bound} also give
\begin{equation}
 |\Delta_t(u)|\leq C_\Delta(M):=44M^2
 \qquad(u\in[0,1]),
 \label{eq:uniform-discriminant-upper-bound}
\end{equation}
for a constant independent of \(Q,t,\delta\).

Completing the square gives
\begin{equation}
 q_u(v)
 =
 c\left(v+\frac{B(u)}{2c}\right)^2
 -\frac{\Delta_t(u)}{4c}.
 \label{eq:uniform-completed-square}
\end{equation}
Let
\[
 V_t(u):=\{v\in[0,1]:|q_u(v)|\leq\delta\}.
\]
We need the following consequences of
\eqref{eq:uniform-completed-square}.

First, if \(\Delta_t(u)<-4|c|\delta\), then \(V_t(u)=\varnothing\).
Indeed, the quadratic has constant sign and its absolute value at its
vertex is \(|\Delta_t(u)|/(4|c|)>\delta\).

Second, if
\begin{equation}
 |\Delta_t(u)|\leq8|c|\delta,
 \label{eq:small-discriminant-condition}
\end{equation}
then
\begin{equation}
 \mu(V_t(u))
 \lesssim_{\alpha,\mu}
 \left(\frac{\delta}{|c|}\right)^{\alpha/2}.
 \label{eq:small-discriminant-fibre}
\end{equation}
To see this, set \(v_0(u):=-B(u)/(2c)\).  From
\eqref{eq:uniform-completed-square},
\[
 |c||v-v_0(u)|^2
 \leq |q_u(v)|+\frac{|\Delta_t(u)|}{4|c|}
 \leq3\delta
\]
whenever \(v\in V_t(u)\).  Thus \(V_t(u)\) is contained in one
interval of length at most \(2\sqrt{3\delta/|c|}\), and
\eqref{eq:small-discriminant-fibre} follows from
\eqref{eq:frostman-uniform}.

Third, if \(\Delta_t(u)>8|c|\delta\), then \(V_t(u)\) is a union of at
most two intervals, each of length at most
\begin{equation}
 \frac{4\delta}{\sqrt{\Delta_t(u)}}.
 \label{eq:large-discriminant-interval-length}
\end{equation}
Indeed, put
\[
 s:=\frac{\sqrt{\Delta_t(u)}}{2|c|},
 \qquad
 \varepsilon:=\frac{\delta}{|c|}.
\]
Then \(s^2>2\varepsilon\), and the sublevel condition implies
\[
 \bigl||v-v_0(u)|^2-s^2\bigr|\leq\varepsilon.
\]
For either possible sign of \(v-v_0(u)\), the admissible radial interval
has length
\begin{align*}
 \sqrt{s^2+\varepsilon}-\sqrt{s^2-\varepsilon}
 &=
 \frac{2\varepsilon}
 {\sqrt{s^2+\varepsilon}+\sqrt{s^2-\varepsilon}}\\
 &\leq\frac{2\varepsilon}{s}
 =\frac{4\delta}{\sqrt{\Delta_t(u)}}.
\end{align*}
This proves \eqref{eq:large-discriminant-interval-length}, and hence
\begin{equation}
 \mu(V_t(u))
 \lesssim_{\alpha,\mu}
 \left(\frac{\delta}{\sqrt{\Delta_t(u)}}\right)^\alpha.
 \label{eq:large-discriminant-fibre}
\end{equation}

We now integrate these fibre estimates.  Define
\[
 U_{-1}:=\{u\in[0,1]:|\Delta_t(u)|\leq8|c|\delta\}
\]
and, for integers \(m\geq0\),
\begin{equation}
 U_m:=\{u\in[0,1]:
          2^{m+3}|c|\delta<\Delta_t(u)
          \leq2^{m+4}|c|\delta\}.
 \label{eq:discriminant-shells}
\end{equation}
Every \(u\) for which \(V_t(u)\neq\varnothing\) belongs to
\(U_{-1}\) or to one of the sets \(U_m\).  By Claim~1, applied to the
quadratic polynomial \(\Delta_t\), whose leading coefficient has
absolute value at least \(\kappa\),
\begin{equation}
 \mu(U_{-1})
 \lesssim_{\alpha,\mu,\kappa}(|c|\delta)^{\alpha/2},
 \qquad
 \mu(U_m)
 \lesssim_{\alpha,\mu,\kappa}
          (2^m|c|\delta)^{\alpha/2}.
 \label{eq:discriminant-shell-masses}
\end{equation}
Combining the first estimate in
\eqref{eq:discriminant-shell-masses} with
\eqref{eq:small-discriminant-fibre} gives
\begin{equation}
 \int_{U_{-1}}\mu(V_t(u))\,d\mu(u)
 \lesssim_{\alpha,\mu,\kappa}\delta^\alpha.
 \label{eq:small-discriminant-contribution}
\end{equation}
For \(u\in U_m\), equation
\eqref{eq:large-discriminant-fibre} gives
\[
 \mu(V_t(u))
 \lesssim_{\alpha,\mu}
 \delta^{\alpha/2}{|c|}^{-\alpha/2}2^{-m\alpha/2}.
\]
Using the second estimate in
\eqref{eq:discriminant-shell-masses}, we therefore obtain
\begin{equation}
 \int_{U_m}\mu(V_t(u))\,d\mu(u)
 \lesssim_{\alpha,\mu,\kappa}
 \delta^\alpha
 \qquad(m\geq0).
 \label{eq:one-discriminant-shell-contribution}
\end{equation}

It remains only to count the nonempty shells.  If \(U_m\neq\varnothing\),
then \eqref{eq:uniform-discriminant-upper-bound} and
\eqref{eq:discriminant-shells} imply
\[
 2^{m+3}|c|\delta<C_\Delta(M).
\]
Since \(|c|\geq\delta\),
\[
 m\leq\log_2(C_\Delta(M))+2\log_2(1/\delta)-3.
\]
Thus only \(O_M(\log(2/\delta))\) shells are nonempty.  Summing
\eqref{eq:small-discriminant-contribution} and
\eqref{eq:one-discriminant-shell-contribution}, and using Fubini's
theorem, gives
\begin{equation}
 (\mu\times\mu)
 \{(u,v):|q_u(v)|\leq\delta\}
 \lesssim_{\alpha,\mu,M,\kappa}
 \delta^\alpha\log(2/\delta).
 \label{eq:large-c-final}
\end{equation}

Equations \eqref{eq:small-c-final} and \eqref{eq:large-c-final}
establish \eqref{eq:uniform-rank2-conclusion} at all sufficiently small
scales, because
\(\log(2/\delta)\leq2\log(1/\delta)\) for
\(0<\delta\leq1/2\). This completes the proof.
\end{proof}

We shall also use the following linear algebra alternative and the
corresponding sublevel estimate.  These two lemmas are needed only to
justify the coordinate-slice step in the proof of
Lemma~\ref{lem:I0_bound}.

\begin{lemma}\label{lem:principal_minor_alternative}
Let
\[
 H
 =
 \begin{pmatrix}
 2d&a&b\\
 a&2e&c\\
 b&c&2g
 \end{pmatrix}
\]
be the Hessian of the quadratic part
\[
 q(x,y,z)
 =
 axy+bxz+cyz+dx^2+ey^2+gz^2.
\]
Assume that $\operatorname{rank}(H)\geq2$.  Then one of the following
alternatives holds.
\begin{enumerate}[label=(\roman*)]
\item
At least one of the coordinate principal minors
\begin{equation}\label{eq:coordinate_principal_minors}
 4de-a^2,
 \qquad
 4dg-b^2,
 \qquad
 4eg-c^2
\end{equation}
is nonzero.  In this case, after permuting the variables, one may
assume that
\[
 4eg-c^2\neq0.
\]
Equivalently, the coordinate slice
\[
 (y,z)\longmapsto q(x,y,z)
\]
has rank-two quadratic part.

\item
All three quantities in
\eqref{eq:coordinate_principal_minors} vanish.  Then necessarily
\begin{equation}\label{eq:sign_incompatible_remaining_case}
 4de=a^2,
 \qquad
 4dg=b^2,
 \qquad
 4eg=c^2,
 \qquad
 abc=-8deg.
\end{equation}
In this case, $H$ has rank three, while every coordinate principal
$2\times2$ minor is zero.  This is the sign-incompatible full-rank
case.
\end{enumerate}
\end{lemma}

\begin{proof}
If one of the three quantities in
\eqref{eq:coordinate_principal_minors} is nonzero, then a permutation
of the variables gives the first alternative.

Assume now that all three quantities vanish:
\begin{equation}\label{eq:all_principal_minors_zero}
 4de=a^2,
 \qquad
 4dg=b^2,
 \qquad
 4eg=c^2.
\end{equation}
We first note that $d,e,g$ are all nonzero.  Indeed, if, for example,
$d=0$, then \eqref{eq:all_principal_minors_zero} gives $a=b=0$, and
the remaining $2\times2$ block has determinant
\[
 4eg-c^2=0.
\]
Hence, $\operatorname{rank}(H)\leq1$, contradicting
$\operatorname{rank}(H)\geq2$.  The same argument applies to $e$ and
$g$.

Multiplying the three identities in
\eqref{eq:all_principal_minors_zero} gives
\[
 64d^2e^2g^2
 =
 a^2b^2c^2
 =
 (abc)^2.
\]
Hence,
\[
 (abc-8deg)(abc+8deg)=0,
\]
so either
\[
 abc=8deg
\]
or
\[
 abc=-8deg.
\]

If $abc=8deg$, then
\[
 c=\frac{ab}{2d},
 \qquad
 e=\frac{a^2}{4d},
 \qquad
 g=\frac{b^2}{4d},
\]
and therefore
\[
 q(x,y,z)
 =
 d\left(
 x+\frac{a}{2d}y+\frac{b}{2d}z
 \right)^2.
\]
Thus, $H$ has rank one, again contradicting
$\operatorname{rank}(H)\geq2$.  Therefore, under the present
hypotheses, the only possible branch with all three principal minors
equal to zero is
\[
 abc=-8deg.
\]

Finally, in this branch $H$ is full rank.  Indeed,
\[
 \det H
 =
 8deg+2abc-2dc^2-2eb^2-2ga^2.
\]
Using \eqref{eq:all_principal_minors_zero}, this simplifies to
\[
 \det H
 =
 2abc-16deg
 =
 2(abc-8deg).
\]
If $abc=-8deg$, then
\[
 \det H=-32deg\neq0.
\]
Hence, $\operatorname{rank}(H)=3$, although all coordinate principal
$2\times2$ minors vanish.
\end{proof}

\begin{lemma}\label{lem:sign_incompatible_sublevel}
Let $\mu_1,\mu_2,\mu_3$ be $\alpha$-Frostman probability measures
supported on fixed bounded intervals in $\mathbb R$.  Let
\[
 Q(x,y,z)
 =
 A_2(x+py+qz)^2+A_1yz+A_0,
\]
where
\[
 A_1A_2\neq0
\]
and $A_0,p,q\in\mathbb R$.  Then, for every
$0<\delta\leq1/2$ and every $t\in\mathbb R$,
\begin{equation}\label{eq:sign_incompatible_sublevel_bound}
 (\mu_1\times\mu_2\times\mu_3)
 \bigl(
 \{
 (x,y,z):
 |Q(x,y,z)-t|\leq\delta
 \}
 \bigr)
 \lesssim_{Q,\mu_1,\mu_2,\mu_3}
 \delta^\alpha.
\end{equation}
\end{lemma}

\begin{proof}
All implicit constants in this proof may depend on $Q$ and on the
Frostman constants and fixed support intervals of the measures.  We
first record a product concentration estimate.  Let $\eta$ be the
pushforward of $\mu_2\times\mu_3$ under
\[
 (y,z)\longmapsto A_1yz.
\]
We claim that, for every $0<\varrho\leq1$ and every $s\in\mathbb R$,
\begin{equation}\label{eq:yz_pushforward_concentration}
\begin{aligned}
 \eta(B(s,\varrho))
 &=
 (\mu_2\times\mu_3)
 \bigl(
 \{
 (y,z):
 |A_1yz-s|\leq\varrho
 \}
 \bigr)\\
 &\lesssim_Q
 \varrho^\alpha
 \log\left(\frac{2}{\varrho}\right).
\end{aligned}
\end{equation}
After absorbing $A_1$ into the implicit constant, fix $y$.
\begin{itemize}
\item
If
\[
 |y|\leq c_Q\varrho
\]
for some constant $c_Q$ depending only on $Q$, we use the trivial bound
in the $z$-variable and the $\alpha$-Frostman estimate in the
$y$-variable, obtaining $O_Q(\varrho^\alpha)$.

\item
If
\[
 2^mc_Q\varrho
 <
 |y|
 \leq
 2^{m+1}c_Q\varrho,
 \qquad
 m\geq0,
\]
then the admissible values of $z$ lie in an interval of length
$O_Q(2^{-m})$, and hence have $\mu_3$-measure
$O_Q(2^{-m\alpha})$.  The corresponding set of $y$ has
$\mu_2$-measure at most
\[
 O_Q((2^m\varrho)^\alpha).
\]
Each nonempty annulus therefore contributes $O_Q(\varrho^\alpha)$,
and there are $O_Q(\log(2/\varrho))$ such annuli because $\mu_2$ is
supported on a bounded interval.
\end{itemize}
Combining both parts proves
\eqref{eq:yz_pushforward_concentration}.

Next, we use a one-dimensional quadratic sublevel estimate.  For fixed
$c,s\in\mathbb R$, set
\[
 E_{c,s}
 :=
 \bigl\{
 x:
 |A_2(x+c)^2-s|\leq\delta
 \bigr\}.
\]
Since $\mu_1$ is $\alpha$-Frostman, we have, uniformly in $c$,
\begin{equation}\label{eq:one_dimensional_quadratic_bound}
 \mu_1(E_{c,s})
 \lesssim_{A_2}
 \delta^\alpha
 (\delta+|s|)^{-\alpha/2}.
\end{equation}
To justify \eqref{eq:one_dimensional_quadratic_bound}, divide by $A_2$
and distinguish
\[
 |s|\lesssim_{A_2}\delta
\]
from
\[
 |s|\gtrsim_{A_2}\delta.
\]
In the first case, the sublevel set is contained in at most two
intervals of total length $O_{A_2}(\sqrt\delta)$.  In the second case,
it is empty if the signs are incompatible; otherwise, by the quadratic
formula, it is contained in at most two intervals of total length
\[
 O_{A_2}\left(
 \frac{\delta}{\sqrt{|s|}}
 \right).
\]
Applying the Frostman bound to these intervals proves
\eqref{eq:one_dimensional_quadratic_bound}, uniformly in the
translation parameter $c$.

Fix $t\in\mathbb R$ and put
\[
 s_0:=t-A_0.
\]
For fixed $y,z$, apply
\eqref{eq:one_dimensional_quadratic_bound} with
\[
 c=py+qz,
 \qquad
 s=s_0-A_1yz.
\]
By Fubini's theorem,
\begin{align*}
 &(\mu_1\times\mu_2\times\mu_3)
 \bigl(
 \{
 (x,y,z):
 |Q(x,y,z)-t|\leq\delta
 \}
 \bigr)\\
 &\qquad
 =
 \iint
 \mu_1(E_{c,s})
 \,d\mu_2(y)\,d\mu_3(z)\\
 &\qquad
 \lesssim_Q
 \int
 \delta^\alpha
 (\delta+|s_0-r|)^{-\alpha/2}
 \,d\eta(r).
\end{align*}
It remains to bound the last integral uniformly in $s_0$.

By \eqref{eq:yz_pushforward_concentration}, the region
\[
 |s_0-r|\leq\delta
\]
contributes at most
\[
\begin{aligned}
 \delta^{\alpha/2}\eta(B(s_0,\delta))
 &\lesssim_Q
 \delta^{3\alpha/2}
 \log\left(\frac2\delta\right)\\
 &\lesssim_Q
 \delta^\alpha.
\end{aligned}
\]
For each integer $m\geq0$ with $2^m\delta<1/2$, consider the truncated
dyadic annulus
\[
 2^m\delta
 <
 |s_0-r|
 \leq
 \min\{2^{m+1}\delta,1/2\}.
\]
It is contained in $B(s_0,2^{m+1}\delta)$, and
$2^{m+1}\delta<1$.  Hence
\eqref{eq:yz_pushforward_concentration}, with
$\varrho=2^{m+1}\delta$, bounds its contribution by
\[
\begin{aligned}
 &\lesssim_Q
 \delta^\alpha
 (2^m\delta)^{-\alpha/2}
 (2^m\delta)^\alpha
 \log\left(\frac{2}{2^m\delta}\right)\\
 &=
 \delta^\alpha
 (2^m\delta)^{\alpha/2}
 \log\left(\frac{2}{2^m\delta}\right).
\end{aligned}
\]
Summing over all such $m$ gives $O_Q(\delta^\alpha)$.  Finally, on the
region
\[
 |s_0-r|>\frac12,
\]
the integrand is bounded by $O(\delta^\alpha)$, so this region also
contributes $O(\delta^\alpha)$.  Combining the three estimates proves
\eqref{eq:sign_incompatible_sublevel_bound}.
\end{proof}

We are now ready to bound $I_0$.

\begin{proof}[Proof of Lemma~\ref{lem:I0_bound}]
Set
\[
 \rho:=\delta^\kappa,
 \qquad
 \mathcal G:=\mathcal G_{2\rho}.
\]
Let
\[
 \mu_3:=\mu\times\mu\times\mu.
\]
By \eqref{eq:def_I0} and Fubini's theorem,
\begin{equation}\label{eq:I0_fubini}
\begin{aligned}
 I_0
 &=
 \int
 \mu_3\bigl(
 \mathcal G
 \cap
 f^{-1}([f(u')-2\delta,f(u')+2\delta])
 \bigr)
 \mathbf 1_{\mathcal G}(u')
 \,d\mu_3(u')\\
 &\leq
 \mu_3(\mathcal G)
 \sup_{t\in\mathbb R}
 \mu_3\bigl(
 \mathcal G
 \cap
 \{|f-t|\leq2\delta\}
 \bigr).
\end{aligned}
\end{equation}

We first bound $\mu_3(\mathcal G)$.  By
Lemmas~\ref{lem:K_dimension} and \ref{lem:Gr_tube}, $K$ is either a
point or an affine line, and
\[
 \mathcal G\subset N_{C\rho}(K).
\]
If $K$ is a point, then Lemma~\ref{lem:tube_mass_line} gives
\begin{equation}\label{eq:G_mass_point}
 \mu_3(\mathcal G)
 \lesssim
 \rho^{3\alpha}.
\end{equation}
If $K$ is an affine line, then Lemma~\ref{lem:tube_mass_line} gives
\begin{equation}\label{eq:G_mass_line}
 \mu_3(\mathcal G)
 \lesssim
 \rho^{2\alpha}.
\end{equation}

We next bound the slice mass uniformly in $t$.  Fix $t\in\mathbb R$
and choose $u_0\in K$.  Since $\nabla f(u_0)=0$, one has
\begin{equation}\label{eq:quadratic_centered}
 f(u)
 =
 f(u_0)
 +
 \frac12(u-u_0)^TH(u-u_0)
\end{equation}
for all $u$, where $H$ is the Hessian from
\eqref{eq:quad_form}.  By Lemma~\ref{lem:K_dimension},
\[
 \operatorname{rank}(H)\in\{2,3\}.
\]

By Lemma~\ref{lem:principal_minor_alternative}, there are two cases.

\medskip
\noindent\emph{Case 1: a coordinate rank-two slice exists.}

After permuting coordinates, we may assume that the principal
$(y,z)$-minor is nonzero:
\[
 \det H_{\{2,3\},\{2,3\}}
 \neq
 0.
\]
For each $x\in[0,1]$, set
\[
 Q_x(y,z):=f(x,y,z).
\]
Then
\[
 D^2Q_x
 =
 D^2_{y,z}f(x,y,z)
 =
 H_{\{2,3\},\{2,3\}},
\]
which is independent of $x$ and satisfies
\[
 \det(D^2Q_x)
 =
 \det(H_{\{2,3\},\{2,3\}})
 \neq
 0.
\]
Thus, the rank-two non-degeneracy of $Q_x$ holds uniformly, with a
lower bound depending only on $f$.  Moreover, the remaining
coefficients of $Q_x$ are uniformly bounded in $x$ because $f$ is
fixed.  Therefore, Lemma~\ref{lem:rank2_sublevel} applies to the family
$\{Q_x:x\in[0,1]\}$ with a constant independent of $x$.

This yields
\[
 (\mu\times\mu)
 \bigl(
 \{
 (y,z)\in[0,1]^2:
 |f(x,y,z)-t|\leq2\delta
 \}
 \bigr)
 \lesssim
 \delta^\alpha
 \log\left(\frac2\delta\right),
\]
uniformly in $x$ and $t$.  Integrating in $x$ gives
\begin{equation}\label{eq:slice_mass}
\begin{aligned}
 \mu_3\bigl(
 \mathcal G\cap\{|f-t|\leq2\delta\}
 \bigr)
 &\leq
 \mu_3\bigl(
 \{|f-t|\leq2\delta\}
 \bigr)\\
 &\lesssim
 \delta^\alpha
 \log\left(\frac2\delta\right),
\end{aligned}
\end{equation}
uniformly in $t$.

Insert \eqref{eq:slice_mass} into \eqref{eq:I0_fubini}, using
\eqref{eq:G_mass_line} in the line case and
\eqref{eq:G_mass_point} in the point case.  We obtain
\[
 I_0
 \lesssim
 \rho^{2\alpha}
 \delta^\alpha
 \log\left(\frac2\delta\right),
\]
where the affine-line bound is used.  Since
\[
 \rho=\delta^\kappa,
\]
this gives
\[
 I_0
 \lesssim
 \delta^{\alpha+2\alpha\kappa}
 \log\left(\frac2\delta\right).
\]
Choosing $\delta_0>0$ so small that
\[
 \log\left(\frac2\delta\right)
 \leq
 \delta^{-\alpha\kappa/2}
 \qquad
 (0<\delta\leq\delta_0),
\]
we obtain the desired conclusion
\[
 I_0
 \lesssim
 \delta^{\alpha+3\alpha\kappa/2}.
\]

\medskip
\noindent\emph{Case 2: the sign-incompatible full-rank case.}

It remains to treat the second alternative in
Lemma~\ref{lem:principal_minor_alternative}, namely
\[
 4de=a^2,
 \qquad
 4dg=b^2,
 \qquad
 4eg=c^2,
 \qquad
 abc=-8deg.
\]
In this case, $H$ is invertible, so
Lemma~\ref{lem:K_dimension} implies that $K$ is a single point.
Hence, by Lemma~\ref{lem:Gr_tube} and the point estimate in
Lemma~\ref{lem:tube_mass_line},
\begin{equation}\label{eq:G_mass_sign_incompatible}
 \mu_3(\mathcal G)
 \lesssim
 \rho^{3\alpha}.
\end{equation}

We next prove the corresponding uniform one-level estimate.  Since
$u_0\in K$, equation \eqref{eq:quadratic_centered} holds.  In the
present sign-incompatible case, set
\[
 p:=\frac{a}{2d},
 \qquad
 q:=\frac{b}{2d}.
\]
Using
\[
 e=\frac{a^2}{4d},
 \qquad
 g=\frac{b^2}{4d},
 \qquad
 c=-\frac{ab}{2d},
\]
the quadratic part can be written as
\[
 d\bigl(
 (x-x_0)+p(y-y_0)+q(z-z_0)
 \bigr)^2
 -
 4dpq\,(y-y_0)(z-z_0),
\]
up to the harmless additive constant $f(u_0)$.  Translating the three
copies of $\mu$ by $x_0,y_0,z_0$ preserves the Frostman exponent and
bounded support.  Therefore,
Lemma~\ref{lem:sign_incompatible_sublevel} gives, uniformly in $t$,
\begin{equation}\label{eq:slice_mass_sign_incompatible}
 \mu_3\bigl(
 \{|f-t|\leq2\delta\}
 \bigr)
 \lesssim
 \delta^\alpha.
\end{equation}
Combining \eqref{eq:I0_fubini},
\eqref{eq:G_mass_sign_incompatible}, and
\eqref{eq:slice_mass_sign_incompatible}, we get
\[
 I_0
 \lesssim
 \rho^{3\alpha}\delta^\alpha
 =
 \delta^{\alpha+3\alpha\kappa}.
\]
Since $0<\delta<1$, this is stronger than
\[
 I_0
 \lesssim
 \delta^{\alpha+3\alpha\kappa/2}.
\]
This completes the proof of the lemma.
\end{proof}

\subsubsection{Conclusion for Ahlfors--David regular measures}

Let $\eta_{\mathrm{AD}}>0$ be the smaller of the transverse gain in
Lemma~\ref{lem:AD-transverse-bound} and $3\alpha\kappa/2$.  The
decomposition \eqref{eq:coincidence-six-term-decomposition} and
Lemma~\ref{lem:I0_bound} give
\[
 \mathcal C_\mu(2\delta)
 \lesssim
 \delta^{\alpha+\eta_{\mathrm{AD}}}.
\]
After replacing $2\delta$ by a general sufficiently small scale,
Lemma~\ref{lem:schwartz-reduction} proves
Theorem~\ref{thm-main-1} in the Ahlfors--David regular case.

\subsection{Passage to general Frostman measures}
\label{sec:Frostman}

We now remove the lower regularity assumption.  Namely, we assume that
$\mu$ is an $\alpha$-Frostman measure, rather than an
$\alpha$-Ahlfors--David regular measure.

Throughout this subsection, $\mu$ is an $\alpha$-Frostman probability
measure supported on $[0,1]$; thus
\begin{equation}\label{eq:frostman_bound}
    \mu(B(x,s))
    \leq
    C_\mu s^\alpha
    \qquad
    \text{for all }x\in\mathbb R\text{ and }s>0.
\end{equation}
We retain the choice of $\kappa$ in
\eqref{eq:transverse-kappa-choice}, as well as the quantities
$I_0,I_1,\ldots,I_6$ defined in
\eqref{eq:def_I0}--\eqref{eq:def_Ij}.  The statements of
Lemmas~\ref{lem:K_dimension}--\ref{lem:sign_incompatible_sublevel}
require only that $\mu$ be $\alpha$-Frostman.  Likewise, the estimate
for $I_0$ in Lemma~\ref{lem:I0_bound} uses only the upper Frostman
bound.  It therefore remains to estimate the six transverse terms
$I_1,\ldots,I_6$ in order to complete the proof of
Theorem~\ref{thm-main-1} for general $\alpha$-Frostman measures.

Fix $0<\delta\leq 1/2$, put
\[
    N:=\lfloor\delta^{-1}\rfloor,
\]
and partition $[0,1]$ into the intervals
\[
    I_m=[m\delta,(m+1)\delta)
    \quad\text{for }0\leq m\leq N-1,
    \qquad
    I_N=[N\delta,1].
\]
Write
\[
    \mathcal D_\delta:=\{I_0,\ldots,I_N\}.
\]
For $I=I_m$, set
\begin{equation}\label{eq:frostman-grid-anchors}
    w(I):=\mu(I),
    \qquad
    a_I:=m\delta.
\end{equation}
Thus, unlike arbitrary points selected from the support, the anchors
form a subset of the fixed grid
$\delta\mathbb Z\cap[0,1]$.  Every cell has diameter at most $\delta$,
and cells of zero mass may be discarded.  For every nonnegative
function $F$ on $[0,1]$,
\begin{equation}\label{eq:discretize_1d}
    \int F\,d\mu
    \leq
    \sum_{I\in\mathcal D_\delta}
    w(I)\sup_{x\in I}F(x).
\end{equation}

Since each cell is contained in an interval of radius $\delta$, the
Frostman bound gives
\[
    w(I)\leq C_\mu\delta^\alpha.
\]
Hence every positive-mass cell belongs to exactly one of the families
\[
    \mathcal D_{\delta,k}
    :=
    \left\{
        I\in\mathcal D_\delta:
        C_\mu2^{-(k+1)}\delta^\alpha
        <
        w(I)
        \leq
        C_\mu2^{-k}\delta^\alpha
    \right\},
    \qquad k\geq0.
\]
Set
\begin{equation}\label{eq:def_Mk}
    M_k
    :=
    \{a_I:I\in\mathcal D_{\delta,k}\}
    \subset
    \delta\mathbb Z\cap[0,1].
\end{equation}

\begin{lemma}\label{lem:level_packing_HHP}
For every $k\geq0$ and every interval $J\subset[0,1]$ with
$|J|\geq\delta$,
\begin{equation}\label{eq:level_packing_HHP}
    |M_k\cap J|
    \leq
    C_\alpha 2^k|J|^\alpha\delta^{-\alpha}.
\end{equation}
\end{lemma}

\begin{proof}
If $a_I\in J$, then $I$ is contained in the $\delta$-neighbourhood
$J^{(\delta)}$ of $J$.  Since $|J|\geq\delta$, this is an interval of
length at most
\[
    |J|+2\delta\leq3|J|.
\]
Because the cells are pairwise disjoint,
\[
\begin{aligned}
    |M_k\cap J|\,C_\mu2^{-(k+1)}\delta^\alpha
    &<
    \sum_{\substack{I\in\mathcal D_{\delta,k}\\ a_I\in J}}
    \mu(I) \\
    &\leq
    \mu(J^{(\delta)}) \\
    &\leq
    C_\mu(3|J|)^\alpha.
\end{aligned}
\]
This proves \eqref{eq:level_packing_HHP}.
\end{proof}

The loss $2^k$ in \eqref{eq:level_packing_HHP} will be removed by a
simple colouring argument.

\begin{lemma}\label{lem:partition_NC_HHP}
Suppose that a nonempty $\delta$-separated set $M\subset[0,1]$
satisfies
\begin{equation}\label{eq:NC_K_HHP}
    |M\cap J|
    \leq
    K|J|^\alpha\delta^{-\alpha}
    \qquad
    (|J|\geq\delta),
\end{equation}
where $K\geq1$.  Then $M$ can be partitioned into at most
\[
    \min\{\lceil2K\rceil+1,|M|\}
\]
nonempty sets $M_\ell'$, each satisfying
\begin{equation}\label{eq:NC_1_HHP}
    |M_\ell'\cap J|
    \leq
    3|J|^\alpha\delta^{-\alpha}
    \qquad
    \text{for every interval }J\text{ with }|J|\geq\delta.
\end{equation}
\end{lemma}

\begin{proof}
Write
\[
    M=\{x_1<\cdots<x_Q\}
\]
and put
\[
    L:=\lceil2K\rceil+1.
\]
Partition the indices $1,\ldots,Q$ according to their residue classes
modulo $L$, and discard the empty residue classes.

The points of $M$ lying in a fixed interval $J$ have consecutive
indices.  Let
\[
    m:=|M\cap J|.
\]
Then every retained residue class contains at most
\[
    \frac{m}{L}+1
\]
of these points.  By \eqref{eq:NC_K_HHP}, and since
$|J|^\alpha\delta^{-\alpha}\geq1$, we have
\[
\begin{aligned}
    \frac{m}{L}+1
    &\leq
    \frac{K|J|^\alpha\delta^{-\alpha}}
         {\lceil2K\rceil+1}
    +1 \\
    &\leq
    \frac12|J|^\alpha\delta^{-\alpha}+1 \\
    &\leq
    \frac32|J|^\alpha\delta^{-\alpha}.
\end{aligned}
\]
This is stronger than \eqref{eq:NC_1_HHP}.  Moreover, the number of
nonempty residue classes is at most
\[
    \min\{L,|M|\},
\]
which proves the required bound.
\end{proof}

Combining Lemmas~\ref{lem:level_packing_HHP} and
\ref{lem:partition_NC_HHP} gives the following formulation, including
the empty-level convention that will be used in the summation.

\begin{corollary}\label{cor:NC_prime_HHP}
Let $\{M_k\}$ be the family of sets defined by \eqref{eq:def_Mk}.
For every $k\geq0$, there is a partition into nonempty classes
\begin{equation}\label{eq:frostman-level-partition}
    M_k
    =
    \bigsqcup_{\ell=1}^{L_k}M_{k,\ell},
\end{equation}
where $L_k=0$ if $M_k=\varnothing$.  Moreover, the partition satisfies
the following two conditions:
\begin{enumerate}
\item
Every $M_{k,\ell}$ is $\delta$-separated and satisfies
\eqref{eq:NC_1_HHP}; in particular,
\[
    |M_{k,\ell}|
    \leq
    3\delta^{-\alpha}.
\]

\item
\begin{equation}\label{eq:frostman-number-of-colours}
    L_k
    \lesssim_\alpha
    \min\{2^k,\delta^{-1}\}.
\end{equation}
\end{enumerate}
\end{corollary}

\begin{proof}
Apply Lemma~\ref{lem:partition_NC_HHP} with
\[
    K=C_\alpha2^k.
\]
After the empty residue classes have been discarded, the number of
classes is also at most
\[
    |M_k|
    \leq
    N+1
    \lesssim
    \delta^{-1}.
\]
Finally, \eqref{eq:NC_1_HHP} with $J=[0,1]$ gives the asserted
cardinality bound.
\end{proof}

\begin{lemma}\label{lem:dyadic_sum_compensation_HHP}
With the convention in Corollary~\ref{cor:NC_prime_HHP},
\begin{equation}\label{eq:dyadic-sum-compensation}
    \sum_{k\geq0}2^{-k}L_k
    \lesssim_\alpha
    \log\left(\frac2\delta\right).
\end{equation}
\end{lemma}

\begin{proof}
Let
\[
    K_\delta
    :=
    \left\lceil
        \log_2\left(\frac1\delta\right)
    \right\rceil.
\]
The first bound in \eqref{eq:frostman-number-of-colours} gives
\[
    \sum_{0\leq k\leq K_\delta}
    2^{-k}L_k
    \lesssim
    K_\delta+1,
\]
whereas the second gives
\[
\begin{aligned}
    \sum_{k>K_\delta}2^{-k}L_k
    &\lesssim
    \delta^{-1}
    \sum_{k>K_\delta}2^{-k} \\
    &\lesssim
    1.
\end{aligned}
\]
Together, these estimates prove
\eqref{eq:dyadic-sum-compensation}.
\end{proof}

\begin{lemma}\label{lem:frostman-transverse-bound}
Let $\eta_{\mathrm{tr}}>0$ be a gain furnished by
Proposition~\ref{prop:six-set-transverse} for the fixed choice of
$\kappa$.  For every $j\in\{1,\ldots,6\}$ and all sufficiently small
$\delta$,
\begin{equation}\label{eq:frostman-all-transverse}
    I_j
    \lesssim
    \delta^{\alpha+\eta_{\mathrm{tr}}/2},
\end{equation}
where the implicit constant is uniform in $j$.
\end{lemma}

\begin{proof}
Recall that the quantities $I_j$ are defined by
\[
    I_j
    :=
    \iint
    \mathbf 1_{\{|f(u)-f(v)|\leq2\delta\}}
    \mathbf 1_{\{|D_j(u,v)|>r\}}
    \,d\mu_3(u)\,d\mu_3(v),
    \qquad
    1\leq j\leq6.
\]
Fix $j\in\{1,\ldots,6\}$.  Apply the cell decomposition in each of the
six variables in the definition of $I_j$.  If a six-tuple of
positive-mass cells contributes, choose a point in the corresponding
cell product at which the two indicator conditions in
\eqref{eq:def_Ij} hold.  Replace its six coordinates by the grid
anchors from \eqref{eq:frostman-grid-anchors}.  As recorded after
\eqref{eq:coincidence-six-term-decomposition}, this replacement
enlarges the coincidence window only to
\begin{equation}\label{eq:frostman-anchor-window}
    |f(a_1,a_2,a_3)-f(a_4,a_5,a_6)|
    \leq
    C_f'\delta
\end{equation}
and leaves the appropriate derivative at the anchor with magnitude at
least
\[
    \frac{r}{2}
    =
    \frac{\delta^\kappa}{2}.
\]

Decompose the six cells into mass levels
$k_1,\ldots,k_6$, and then use the partitions
\eqref{eq:frostman-level-partition}.  For fixed
$(k_q,\ell_q)_{q=1}^6$, assign the six sets
$M_{k_q,\ell_q}$ to the six coordinates in their original order.
These sets need not be equal and need not have equal cardinalities.
Empty levels contribute nothing.  By
Corollary~\ref{cor:NC_prime_HHP}, every nonempty set is
$\delta$-separated and satisfies
\[
    |M_{k_q,\ell_q}\cap J|
    \leq
    3|J|^\alpha\delta^{-\alpha}
    \qquad
    \text{for every interval }J\text{ with }|J|\geq\delta.
\]

Let
\[
    \Sigma_j(k_1,\ell_1;\ldots;k_6,\ell_6)
\]
be the number of anchor tuples satisfying
\eqref{eq:frostman-anchor-window} and the relevant anchored derivative
condition.  Proposition~\ref{prop:six-set-transverse}, with
\[
    C_0=3,
    \qquad
    C_1=C_f',
    \qquad
    c_1=\frac12,
\]
applies directly to these six possibly unequal sets.  Its final
assertion covers $j=4,5,6$, when the derivative is imposed on the
unprimed triple.  Thus, uniformly in all levels and classes,
\begin{equation}\label{eq:Sigma_restricted_bound_HHP}
    \Sigma_j(k_1,\ell_1;\ldots;k_6,\ell_6)
    \lesssim
    \delta^{-5\alpha+\eta_{\mathrm{tr}}}.
\end{equation}

By \eqref{eq:def_Mk}, a cell in level $k_q$ has mass at most
\[
    C_\mu2^{-k_q}\delta^\alpha.
\]
Summing \eqref{eq:Sigma_restricted_bound_HHP} over the six levels and
their colour classes therefore gives
\begin{align*}
    I_j
    &\lesssim
    C_\mu^6\delta^{6\alpha}
    \sum_{k_1,\ldots,k_6\geq0}
    2^{-(k_1+\cdots+k_6)}
    \sum_{\substack{
        1\leq\ell_1\leq L_{k_1}\\
        \cdots\\
        1\leq\ell_6\leq L_{k_6}
    }}
    \delta^{-5\alpha+\eta_{\mathrm{tr}}}
    \\
    &\lesssim
    C_\mu^6\delta^{\alpha+\eta_{\mathrm{tr}}}
    \left(
        \sum_{k\geq0}2^{-k}L_k
    \right)^6
    \\
    &\lesssim
    C_\mu^6\delta^{\alpha+\eta_{\mathrm{tr}}}
    \log^6\left(\frac2\delta\right)
    \\
    &\lesssim
    \delta^{\alpha+\eta_{\mathrm{tr}}/2},
\end{align*}
after decreasing $\delta_0$ if necessary.  Since the argument applies
to every $j$, this proves all six estimates in
\eqref{eq:frostman-all-transverse}.
\end{proof}

We can now complete the proof of the main theorem.

\begin{proof}[Proof of Theorem~\ref{thm-main-1}]
Set
\begin{equation}\label{eq:final-coincidence-gain}
    \eta
    :=
    \min\left\{
        \frac{\eta_{\mathrm{tr}}}{2},
        \frac{3\alpha\kappa}{2}
    \right\}
    >0.
\end{equation}
Lemma~\ref{lem:I0_bound},
Lemma~\ref{lem:frostman-transverse-bound}, and the decomposition
\eqref{eq:coincidence-six-term-decomposition} give
\begin{equation}\label{eq:frostman-coincidence-at-grid-scale}
    \mathcal C_\mu(2\delta)
    \lesssim
    \delta^{\alpha+\eta}
\end{equation}
for all sufficiently small $\delta$.  Given
$0<\rho\leq\rho_0$, with $\rho_0$ sufficiently small, take
$\delta=\rho/2$ in
\eqref{eq:frostman-coincidence-at-grid-scale}.  We obtain precisely
\eqref{eq:coincidence-target}, namely,
\[
    \mathcal C_\mu(\rho)
    \lesssim
    \rho^{\alpha+\eta}.
\]
Finally, Lemma~\ref{lem:schwartz-reduction}, applied with this
coincidence bound, yields
\[
    \left\|
        \varphi_\delta*f_\#(\mu^{\otimes3})
    \right\|_{L^2(\mathbb R)}^2
    \lesssim
    \delta^{\alpha+\eta-1}.
\]
This proves Theorem~\ref{thm-main-1} with
$\varepsilon=\eta$.
\end{proof}

\section{An application to level-set dimensions}
\label{subsec:level-set-dimensions}
We record an application of Theorem~\ref{thm-main-1} to the size of level
sets of the quadratic image map.  The argument is a regular-value-type
principle: a quantitative upper bound for the mass of a thin neighbourhood
of a level set gives an upper box dimension bound for that level set.

For a bounded set \(E\subset\mathbb R^d\), let \(\mathscr C_r(E)\) denote the
smallest number of Euclidean balls of radius \(r\) required to cover \(E\), and
define its upper box dimension by
\[
    \overline{\dim}_{\mathrm B}E
    :=
    \limsup_{r\downarrow0}
    \frac{\log \mathscr C_r(E)}{-\log r}.
\]
We use the convention
\[
    \overline{\dim}_{\mathrm B}\emptyset=-\infty.
\]
If \(A\subset[0,1]\), \(\mu\) is a probability measure with
\(\spt(\mu)=A\), and \(f\in\mathbb R[x,y,z]\), set
\[
    \nu:=f_\#(\mu\times\mu\times\mu),
\]
and, for \(t\in\mathbb R\), define the level set
\[
    \mathcal{F}_t:=\{(x,y,z)\in A^3:f(x,y,z)=t\}.
\]
\begin{theorem}
\label{thm:level-set-typical-exceptional}
Let \(f\in\mathbb R[x,y,z]\) be a fixed non-degenerate quadratic polynomial, and let \(\mu\) be an \(\alpha\)-AD-regular
probability measure on \(A=\spt(\mu)\subset[0,1]\), with
\(0<\alpha<1\).  Let \(\varphi\) be a fixed nonnegative Schwartz
function satisfying \(\int_{\mathbb R}\varphi=1\), and put
\(\varphi_\delta(t)=\delta^{-1}\varphi(t/\delta)\). Assume that, for some \(\sigma\in(\alpha,1]\),
\begin{equation}\label{eq:level-set-L2-hypothesis}
    \|\varphi_\delta*\nu\|_{L^2(\mathbb R)}^2
    \le C_\sigma \delta^{\sigma-1},
\end{equation}
for all sufficiently small \(\delta>0\).  Then the following hold.
\begin{enumerate}[label=\textup{(\roman*)}]
    \item For \(\nu\)-almost every \(t\in\mathbb R\),
    \[
        \overline{\dim}_{\mathrm B}\mathcal{F}_t\le 3\alpha-\sigma.
    \]

    \item For every \(\beta\in(0,1)\), set
    \[
        \mathcal E_\beta
        :=
        \{t\in\mathbb R:
        \overline{\dim}_{\mathrm B}\mathcal{F}_t>3\alpha-\beta\}.
    \]
    If \(0<\beta\le\alpha\), then \(\mathcal E_\beta=\emptyset\).  If
    \(\alpha<\beta<1\), then
    \[
        \dim_{\mathrm H}\mathcal E_\beta
        \le
        \max\{0,2\beta-\sigma\}.
    \]
\end{enumerate}
In particular, if \(\varepsilon>0\) is the gain supplied by
Theorem~\ref{thm-main-1}, then \eqref{eq:level-set-L2-hypothesis} holds for
every \(\sigma\) with
\[
    \alpha<\sigma\le \min\{1,\alpha+\varepsilon\}.
\]
\end{theorem}
We first record the following packing lemma: 
\begin{lemma}
\label{lem:level-set-packing-local-mass}
Let \(\mu\) be \(\alpha\)-AD-regular on \(A=\spt(\mu)\subset[0,1]\), and let
\(f:[0,1]^3\to\mathbb R\) be Lipschitz. Suppose that, for some \(t\in\mathbb R\) and some \(\theta>0\), there are
constants \(C_t>0\) and \(r_t>0\) such that
\begin{equation}\label{eq:local-value-distribution-level-set}
    \nu(B(t,r))\le C_t r^\theta
    \qquad\text{for all }0<r<r_t.
\end{equation}
Then
\[
    \overline{\dim}_{\mathrm B}\mathcal{F}_t\le 3\alpha-\theta.
\]
\end{lemma}

\begin{proof}
If \(\mathcal{F}_t=\emptyset\), there is nothing to prove.  Let \(L\ge1\) be a Lipschitz
constant for \(f\) on \([0,1]^3\).  For \(r>0\) sufficiently small, choose a maximal $3r$-separated subset
\[
    \{p_1,\ldots,p_n\}\subset \mathcal{F}_t.
\]
Then the balls $B(p_j,r)$ are pairwise disjoint, and maximality
implies that the balls $B(p_j,3r)$ cover $\mathcal{F}_t$.  Thus
\[\mathscr C_{3r}(\mathcal{F}_t)\le n.\]

Since \(\mu\) is \(\alpha\)-AD-regular on \(A\), the product measure
\(\mu\times\mu\times\mu\) satisfies
\begin{equation}\label{eq:product-lower-ad-level-set}
    (\mu\times\mu\times\mu)(B(p,r))
    \ge c r^{3\alpha},
\end{equation}
for every \(p\in A^3\) and all sufficiently small \(r>0\).  Indeed, if
\(p=(x,y,z)\), then
\[B(x,r/\sqrt3)\times B(y,r/\sqrt3)\times B(z,r/\sqrt3)\subset B(p,r),\]
and the lower AD-regularity of \(\mu\) gives
\eqref{eq:product-lower-ad-level-set}.

For each \(q\in B(p_j,r)\cap A^3\), since \(f(p_j)=t\), we have
\[|f(q)-t|\le Lr.\]
Therefore, by disjointness, \eqref{eq:product-lower-ad-level-set}, and
\eqref{eq:local-value-distribution-level-set}, provided \(r\) is small enough
that \(Lr<r_t\),
\[ncr^{3\alpha}\le(\mu\times\mu\times\mu)\{q\in A^3:|f(q)-t|\le Lr\}=\nu(B(t,Lr))\le C_t L^\theta r^\theta.\]
It follows that
\[\mathscr C_{3r}(\mathcal{F}_t)\le c^{-1}C_tL^\theta r^{-(3\alpha-\theta)},\]
for all sufficiently small \(r>0\).  Taking the upper limit as \(r\downarrow0\) proves the lemma.
\end{proof}

\begin{proposition}
\label{prop:level-set-uniform-bound}
Let \(f\in\mathbb R[x,y,z]\) be a fixed non-degenerate quadratic polynomial, and let \(\mu\) be an \(\alpha\)-AD-regular
probability measure on \(A=\spt(\mu)\subset[0,1]\), with
\(0<\alpha<1\).  For \(t\in\mathbb R\), set
\[
    \mathcal{F}_t:=\{(x,y,z)\in A^3:f(x,y,z)=t\}.
\]
Then, for every \(t\in\mathbb R\),
\[
    \overline{\dim}_{\mathrm B}\mathcal{F}_t\le 2\alpha.
\]
\end{proposition}

\begin{remark}[Sharpness of the uniform level-set bound]
The exponent \(2\alpha\) in
Proposition~\ref{prop:level-set-uniform-bound} is sharp in general.
Indeed, consider the quadratic polynomial
\[
    f(x,y,z):=(x-y)(x-z)
             =x^2-xy-xz+yz.
\]
This polynomial depends nontrivially on all three variables and is
non-degenerate.

For \(t=0\), the corresponding level set is
\[
\begin{aligned}
    \mathcal F_0
    &=\{(x,y,z)\in A^3:(x-y)(x-z)=0\}\\
    &=\{(x,x,z):x,z\in A\}
      \,\cup\,
      \{(x,y,x):x,y\in A\}.
\end{aligned}
\]
Each of the two sets on the right is bi-Lipschitz equivalent to
\(A\times A\).  Since \(A\) is \(\alpha\)-AD-regular, \(A\times A\) is
\(2\alpha\)-AD-regular and hence has upper box dimension \(2\alpha\).
Consequently,
\[
    \overline{\dim}_{\mathrm B}\mathcal F_0=2\alpha.
\]
Therefore, the exponent in
Proposition~\ref{prop:level-set-uniform-bound} cannot be improved uniformly
over all non-degenerate quadratic polynomials.
\end{remark}
\begin{proof}[Proof of Proposition \ref{prop:level-set-uniform-bound}]
Set
\[
    \nu:=f_\#(\mu\times\mu\times\mu).
\]
We first recall the one-level concentration estimate
\begin{equation}\label{eq:level-set-uniform-one-level-concentration}
    \nu(B(t,r))
    \lesssim
    r^\alpha\log\Bigl(\frac{2}{r}\Bigr),
    \qquad 0<r\le\frac12,
\end{equation}
uniformly in \(t\in\mathbb R\).  We include the proof below, since this is the only input needed for the uniform level-set estimate.

Let \(u=(x,y,z)^T\), and write
\[f(u)=\tfrac12u^T H u+\vec{b}\cdot u+c_0,\qquad\nabla f(u)=Hu+\vec{b},
\]
as in \eqref{eq:quad_form}, and put
\[\mathcal{K}:=\{u\in\mathbb R^3:\nabla f(u)=0\}.\]

First suppose that \(\mathcal{K}=\emptyset\).  Since \([0,1]^3\) is compact, there exists
\(c_f>0\) such that \(|\nabla f(u)|\ge c_f\) on \([0,1]^3\).  The cube is
covered by the three sets on which one of $|\partial_x f|,|\partial_y f|,|\partial_z f|$ is at least \(c_f/\sqrt3\).  Consider, for instance, the set where
\(|\partial_x f|\ge c_f/\sqrt3\).  For fixed \((y,z)\), the one-variable
polynomial \(x\mapsto f(x,y,z)\) is quadratic, and on the set where
\(|\partial_x f|\ge c_f/\sqrt3\), it is monotone on \(O_f(1)\) intervals with
slope bounded from below.  Hence
\[
    \{x\in[0,1]: |f(x,y,z)-t|\le r,\ 
    |\partial_x f(x,y,z)|\ge c_f/\sqrt3\}
\]
is contained in a union of \(O_f(1)\) intervals of length \(O_f(r)\).  The
Frostman condition gives that its \(\mu\)-measure is \(O(r^\alpha)\), uniformly
in \(y,z,t\).  Integrating in \((y,z)\) and summing over the three coordinate
choices gives
\[
    \nu(B(t,r))\lesssim r^\alpha.
\]

It remains to consider the case \(\mathcal{K}\ne\emptyset\).  By
Lemma~\ref{lem:K_dimension}, \(\operatorname{rank}(H)\in\{2,3\}\).  We apply
Lemma~\ref{lem:principal_minor_alternative} and consider its two alternatives:
\begin{itemize}
\item If a coordinate rank-two slice
exists, then, after permuting the variables, we may assume that the principal
\((y,z)\)-minor of \(H\) is invertible.  For fixed \(x\in[0,1]\), set
\[
    Q_x(y,z):=f(x,y,z).
\]
Then
\[
    D^2_{y,z}Q_x=H_{\{2,3\},\{2,3\}},
\]
which is independent of \(x\) and has rank two.  Applying
Lemma~\ref{lem:rank2_sublevel} to \(Q_x\), with constants uniform in \(x\), gives
\[
    (\mu\times\mu)\{(y,z)\in[0,1]^2: |f(x,y,z)-t|\le r\}
    \lesssim
    r^\alpha\log\Bigl(\frac{2}{r}\Bigr),
\]
uniformly in \(x\) and \(t\).  Integrating in \(x\) gives
\eqref{eq:level-set-uniform-one-level-concentration} in this case.
\item Assume that all three coordinate principal
$2\times2$ minors vanish. In the notation
\[
    q(x,y,z)=axy+bxz+cyz+dx^2+ey^2+gz^2
\]
for the quadratic part of \(f\), this case is
\[
    4de=a^2,\qquad
    4dg=b^2,\qquad
    4eg=c^2,\qquad
    abc=-8deg.
\]
The Hessian is invertible, so \(\mathcal{K}\) consists of a single point
\(u_0=(x_0,y_0,z_0)\).  Since \(\nabla f(u_0)=0\),
\[
    f(u)=f(u_0)+\tfrac12(u-u_0)^T H(u-u_0).
\]
Setting
\[
    p_1:=\frac{a}{2d},
    \qquad
    p_2:=\frac{b}{2d},
\]
and using
\[
    e=\frac{a^2}{4d},\qquad
    g=\frac{b^2}{4d},\qquad
    c=-\frac{ab}{2d},
\]
the quadratic part in centred variables has the form
\[
    d\bigl((x-x_0)+p_1(y-y_0)+p_2(z-z_0)\bigr)^2
    -
    4dp_1p_2\,(y-y_0)(z-z_0).
\]
After translating the three copies of \(\mu\) by \(x_0,y_0,z_0\), the measures
remain \(\alpha\)-Frostman and supported on fixed bounded intervals.  Therefore
Lemma~\ref{lem:sign_incompatible_sublevel} gives the stronger estimate
\[
    \nu(B(t,r))\lesssim r^\alpha,
\]
uniformly in \(t\).  This completes the proof of
\eqref{eq:level-set-uniform-one-level-concentration}.
\end{itemize}
Now fix \(0<\eta<\alpha\).  Since
\[
    \log\Bigl(\frac{2}{r}\Bigr)\le C_\eta r^{-\eta},
\]
for all sufficiently small \(r>0\), \eqref{eq:level-set-uniform-one-level-concentration}
implies
\[
    \nu(B(t,r))\le C_\eta r^{\alpha-\eta},
\]
uniformly in \(t\).  Applying Lemma~\ref{lem:level-set-packing-local-mass} with
\(\theta=\alpha-\eta\), we obtain
\[
    \overline{\dim}_{\mathrm B}\mathcal{F}_t
    \le
    3\alpha-(\alpha-\eta)
    =
    2\alpha+\eta.
\]
Letting \(\eta\downarrow0\) proves the proposition. 
\end{proof}

\begin{proof}[Proof of Theorem \ref{thm:level-set-typical-exceptional}]
We first derive a correlation estimate from \eqref{eq:level-set-L2-hypothesis}.
Let
\[
    \widetilde\varphi(u):=\varphi(-u),
    \qquad
    K_0:=\varphi*\widetilde\varphi.
\]
Since \(\varphi\ge0\) and \(\int\varphi=1\), one has \(K_0\ge0\) and
\[
    K_0(0)=\int\varphi(u)^2\,du>0.
\]
Thus there are constants \(c_0>0\) and \(a_0>0\) such that
\[
    K_0(u)\ge c_0
    \qquad\text{for all } |u|\le a_0.
\]
If \(K_\delta:=\varphi_\delta*\widetilde\varphi_\delta\), then
\[
    K_\delta(u)
    =
    \delta^{-1}K_0(u/\delta)
    \ge
    c_0\delta^{-1}\mathbf 1_{\{|u|\le a_0\delta\}}.
\]
By Fubini's theorem,
\[
    \|\varphi_\delta*\nu\|_{L^2(\mathbb R)}^2
    =
    \iint K_\delta(u-v)\,d\nu(u)d\nu(v).
\]
Therefore, 
\[
    (\nu\times\nu)\{(u,v):|u-v|\le a_0\delta\}
    \lesssim\delta\|\varphi_\delta*\nu\|_{L^2(\mathbb R)}^2
    \le C_\sigma\delta^\sigma.
\]
Applying the preceding estimate with \(\delta=r/a_0\) shows that, after
altering the constant and decreasing the admissible upper scale,
\begin{equation}\label{eq:level-set-correlation-bound}
    (\nu\times\nu)\{(u,v):|u-v|\le r\}
    \lesssim r^\sigma,
\end{equation}
for every sufficiently small \(r>0\). 

Letting $r\downarrow0$ in
\eqref{eq:level-set-correlation-bound} and using continuity from above gives
\[
    (\nu\times\nu)(\{(u,u):u\in\mathbb R\})=0.
\]
Thus the diagonal makes no contribution to the layer-cake formula.

Let \(0<\gamma<\sigma\).  Since \(f([0,1]^3)\) is compact,
\eqref{eq:level-set-correlation-bound} and the layer-cake representation for
Riesz kernels imply
\begin{align}\label{eq:level-set-finite-energy}
I_\gamma(\nu)&:=\iint |u-v|^{-\gamma}\,d\nu(u)d\nu(v) \notag\\
&\lesssim1+\gamma\int_0^1 r^{-\gamma-1}(\nu\times\nu)\{(u,v):|u-v|\le r\}\,dr \notag\\
&\lesssim1+\gamma\int_0^1 r^{\sigma-\gamma-1}\,dr<\infty.
\end{align}
Thus the potential
\[
    U_\gamma(t):=\int |t-u|^{-\gamma}\,d\nu(u)
\]
is finite for \(\nu\)-almost every \(t\).  For such \(t\),
\[\nu(B(t,r))\le U_\gamma(t)r^\gamma\text{ for all }0<r<1.\]
By Lemma~\ref{lem:level-set-packing-local-mass},
\[
    \overline{\dim}_{\mathrm B}\mathcal{F}_t\le 3\alpha-\gamma
\]
for \(\nu\)-almost every \(t\).  Taking a countable sequence of rational
\(\gamma\uparrow\sigma\) proves (i).

We now prove (ii).  For \(0<\beta\le\alpha\), Proposition~\ref{prop:level-set-uniform-bound}
gives
\[
    \overline{\dim}_{\mathrm B}\mathcal{F}_t\le 2\alpha\le 3\alpha-\beta
    \qquad\text{for every }t\in\mathbb R.
\]
Therefore \(\mathcal E_\beta=\emptyset\).

It remains to consider \(\alpha<\beta<1\).  Since \(A=\spt(\mu)\) and \(f\) is
continuous, if \(\mathcal{F}_t\ne\emptyset\), then \(t\in\spt(\nu)\).  Thus
\[
    \mathcal E_\beta\subset\spt(\nu).
\]
For \(t\in\spt(\nu)\), define
\[
    \underline d_\nu(t)
    :=
    \liminf_{r\downarrow0}
    \frac{\log\nu(B(t,r))}{\log r}.
\]
If \(\underline d_\nu(t)>\beta\), then
\[
    \nu(B(t,r))\le r^\beta
\]
for all sufficiently small \(r>0\), and Lemma~\ref{lem:level-set-packing-local-mass}
gives
\[
    \overline{\dim}_{\mathrm B}\mathcal{F}_t\le 3\alpha-\beta.
\]
Hence
\begin{equation}\label{eq:level-set-bad-contained-localdim}
    \mathcal E_\beta
    \subset
    E_\beta
    :=
    \{t\in\spt(\nu):\underline d_\nu(t)\le\beta\}.
\end{equation}

Fix \(0<\gamma<\sigma\).  By \eqref{eq:level-set-finite-energy},
\(I_\gamma(\nu)<\infty\).  The truncated energy
\[
    I_\gamma(\nu;\rho)
    :=
    \iint_{0<|u-v|\le\rho}|u-v|^{-\gamma}\,d\nu(u)d\nu(v)
\]
tends to \(0\) as \(\rho\downarrow0\) by dominated convergence.  Let
\(\eta>0\) and \(\rho>0\).  Define
\(E_{\beta,\eta,\rho}\) to be the set of all \(t\in E_\beta\) for which there
exists \(0<r_t<\rho\) with
\[
    \nu(B(t,r_t))\ge r_t^{\beta+\eta}.
\]
By the definition of lower local dimension,
\[
    E_\beta\subset E_{\beta,\eta,\rho}
\]
for every \(\eta>0\) and every \(\rho>0\).  Applying the Vitali covering lemma
to the family \(\{B(t,r_t):t\in E_{\beta,\eta,\rho}\}\), we obtain a countable
pairwise disjoint subcollection \(B_i=B(t_i,r_i)\) such that
\[
    E_{\beta,\eta,\rho}\subset\bigcup_i B(t_i,5r_i).
\]
Put
\[
    s=2(\beta+\eta)-\gamma
\]
and assume \(s>0\).  Then
\[
    \mathcal H^s_{10\rho}(E_\beta)
    \le
    \mathcal H^s_{10\rho}(E_{\beta,\eta,\rho})
    \le
    C\sum_i r_i^s
    =
    C\sum_i r_i^{2(\beta+\eta)-\gamma}
    \le
    C\sum_i r_i^{-\gamma}\nu(B_i)^2.
\]
Since the intervals \(B_i\) are pairwise disjoint and \(|u-v|\le2r_i\le2\rho\)
on \(B_i\times B_i\), we also have
\[
    I_\gamma(\nu;2\rho)
    \ge
    \sum_i
    \iint_{B_i\times B_i}|u-v|^{-\gamma}\,d\nu(u)d\nu(v)
    \ge
    \sum_i(2r_i)^{-\gamma}\nu(B_i)^2.
\]
Therefore
\[
    \mathcal H^s_{10\rho}(E_\beta)
    \le
    C I_\gamma(\nu;2\rho).
\]
Letting \(\rho\downarrow0\) gives \(\mathcal H^s(E_\beta)=0\), and hence
\[\dim_{\mathrm H}E_\beta\le 2(\beta+\eta)-\gamma.\]
\begin{itemize}
\item If \(2\beta\ge\sigma\), we let \(\gamma\uparrow\sigma\) and then
\(\eta\downarrow0\) to obtain
\[\dim_{\mathrm H}E_\beta\le 2\beta-\sigma\]
and thus $\dim_{\mathrm H}\mathcal{E}_\beta\leq2\beta-\sigma$.
\item If \(2\beta<\sigma\), we first choose \(\eta>0\) so small that
\(2(\beta+\eta)<\sigma\).  We then choose an arbitrarily small
\(s\in(0,2(\beta+\eta))\) and set
\[\gamma=2(\beta+\eta)-s.\]
This gives
\[\dim_{\mathrm H}E_\beta\le s,\]
hence \(\dim_{\mathrm H}E_\beta=0\).  Combining this with
\eqref{eq:level-set-bad-contained-localdim}, we complete the proof of part (ii).
\end{itemize}

\end{proof}
\section{A concluding remark}\label{section4448}

In this section, we record an example showing that the power-saving estimate in
Theorem~\ref{thm-main-1} does not, in general, imply that the smoothed
\(L^2\)-energy remains bounded as the scale tends to zero.  
\begin{proposition}\label{prop:unboundedness}
Let \(0<\alpha<1/2\), and let \(\mu\) be an \(\alpha\)-AD-regular
probability measure on \(A:=\spt(\mu)\subset[0,1]\).  Define
\[
 f(x,y,z)=x+(y-z)^2,
 \qquad
 \nu=f_\#(\mu\times\mu\times\mu).
\]
For every fixed nonnegative Schwartz function \(\varphi\not\equiv0\), put
\(\varphi_\delta(t):=\delta^{-1}\varphi(t/\delta)\).  Then
\[
 \|\varphi_\delta*\nu\|_{L^2(\mathbb R)}^2
 \gtrsim\delta^{2\alpha-1},
\]
for all sufficiently small \(\delta>0\).  In particular, the left-hand side
tends to infinity as \(\delta\downarrow0\).
\end{proposition}

\begin{proof}
The polynomial $f(x,y,z)=x+(y-z)^2$ is non-degenerate: if
$f=G(I_1(x)+I_2(y)+I_3(z))$, then the identity $f_x=1$ forces $G$ and
$I_1$ to be affine, whereas $f_{yz}=-2$ rules out additive
separability.
Put \(\widetilde\varphi(t)=\varphi(-t)\) and
\(K=\varphi*\widetilde\varphi\).  Then \(K\ge0\) and
\(K(0)=\|\varphi\|_2^2>0\).  Hence there are \(a_0,c_0>0\) such that
\begin{equation}\label{eq:kernel-positive-near-zero-unboundedness}
 K(s)\ge c_0\qquad(|s|\le a_0).
\end{equation}

Choose a fixed \(c\in(0,1)\), and let \(\lambda_\delta\) be the submeasure
of \(\nu\) defined by
\[
 \lambda_\delta(B)
 :=\iiint
 \mathbf 1_{\{|y-z|\le c\sqrt\delta\}}
 \mathbf 1_B\bigl(x+(y-z)^2\bigr)
 \,d\mu(x)d\mu(y)d\mu(z).
\]
The lower Ahlfors-David bound gives
\begin{equation}\label{eq:lambda-total-mass-unboundedness}
 \|\lambda_\delta\|
 =\int\mu(B(z,c\sqrt\delta))\,d\mu(z)
 \gtrsim\delta^{\alpha/2}.
\end{equation}
Moreover, by construction,
\[\spt(\lambda_\delta)\subset A+[0,c^2\delta].\]

Set \(h=a_0/2\), and partition \(\mathbb R\) into half-open intervals of
length \(h\delta\).  We claim that at most \(C\delta^{-\alpha}\) of these
intervals meet \(A+[0,c^2\delta]\).  Indeed, let \(\{x_1,\ldots,x_N\}\)
be a maximal \(h\delta\)-separated subset of \(A\).  The balls
\(B(x_j,h\delta/3)\) are disjoint, and lower Ahlfors-David regularity gives
\[
 1\ge\sum_{j=1}^N\mu(B(x_j,h\delta/3))
 \gtrsim N\delta^\alpha.
\]
Thus \(N\lesssim\delta^{-\alpha}\).  Maximality covers \(A\) by the balls
\(B(x_j,h\delta)\), and enlarging by \([0,c^2\delta]\) increases the number
of length-\(h\delta\) intervals by only a fixed factor.  This proves the
claim.

Let \(\mathcal Q_\delta\) be the intervals in this partition that meet
\(\spt(\lambda_\delta)\).  By the Cauchy-Schwarz inequality and
\eqref{eq:lambda-total-mass-unboundedness},
\begin{equation}\label{eq:lambda-square-sum-unboundedness}
 \sum_{Q\in\mathcal Q_\delta}\lambda_\delta(Q)^2
 \ge\frac{\|\lambda_\delta\|^2}{|\mathcal Q_\delta|}
 \gtrsim\delta^{2\alpha}.
\end{equation}

Since \(\lambda_\delta\le\nu\) and \(K\ge0\), the kernel identity and
\eqref{eq:kernel-positive-near-zero-unboundedness} imply
\begin{align*}
 \|\varphi_\delta*\nu\|_2^2
 &=\delta^{-1}\iint K((u-v)/\delta)\,d\nu(u)d\nu(v)\\
 &\ge\delta^{-1}\iint K((u-v)/\delta)
       \,d\lambda_\delta(u)d\lambda_\delta(v)\\
 &\ge c_0\delta^{-1}
       \sum_{Q\in\mathcal Q_\delta}\lambda_\delta(Q)^2
 \gtrsim\delta^{2\alpha-1},
\end{align*}
because two points in the same \(Q\) are at distance at most
\(h\delta<a_0\delta\).  Since \(2\alpha-1<0\), the final assertion follows.
\end{proof}


\begin{thebibliography}{99}

\footnotesize
\setlength{\itemsep}{0pt}
\setlength{\parskip}{0pt}
\bibitem{chao}
N.~Arala and S.~Chow,
\textit{Expansion properties of polynomials over finite fields},
Finite Fields and Their Applications \textbf{108} (2025),
Article~102687.

\bibitem{benet}
M.~Bennett, D.~Hart, A.~Iosevich, J.~Pakianathan, and M.~Rudnev,
\textit{Group actions and geometric combinatorics in \(\mathbb F_q^d\)},
Forum Mathematicum \textbf{29}(1) (2017), 91--110.

\bibitem{chap}
J.~Chapman, M.~B.~Erdo\u{g}an, D.~Hart, A.~Iosevich, and D.~Koh,
\textit{Pinned distance sets, \(k\)-simplices, Wolff's exponent in finite
fields and sum-product estimates},
Mathematische Zeitschrift \textbf{271}(1--2) (2012), 63--93.

\bibitem{DemeterORegan}
C.~Demeter and W.~O'Regan,
\textit{Incidence estimates for quasi-product sets and applications},
arXiv preprint arXiv:2511.15899v1 (2025).

\bibitem{Eiosevich}
S.~Eswarathasan, A.~Iosevich, and K.~Taylor,
\textit{Fourier integral operators, fractal sets, and the regular value
theorem},
Advances in Mathematics \textbf{228}(4) (2011), 2385--2402.

\bibitem{FuRen}
Y.~Fu and K.~Ren,
\textit{Incidence estimates for \(\alpha\)-dimensional tubes and
\(\beta\)-dimensional balls in \(\mathbb R^2\)},
Journal of Fractal Geometry \textbf{11}(1--2) (2024), 1--30.
\bibitem{gra}
L. Grafakos, A. Greenleaf, A. Iosevich, and E. Palsson, \textit{Multilinear generalized Radon transforms and point configurations}, Forum Mathematicum, \textbf{27}(4) (2015), 2323--2360. 
\bibitem{GIT}
A.~Greenleaf, A.~Iosevich, and K.~Taylor,
\textit{On \(k\)-point configuration sets with nonempty interior},
Mathematika \textbf{68}(1) (2022), 163--190.

\bibitem{ha2}
D.~Hart and A.~Iosevich,
\textit{Sums and products in finite fields: an integral geometric viewpoint},
in \textit{Radon Transforms, Geometry, and Wavelets},
Contemporary Mathematics \textbf{464} (2008), 129--135.

\bibitem{KPS}
D.~Koh, T.~Pham, and C.-Y.~Shen,
\textit{Falconer type functions in three variables},
Journal of Functional Analysis \textbf{286}(4) (2024),
Article~110246.

\bibitem{MP}
B.~Murphy and G.~Petridis,
\textit{Products of differences over arbitrary finite fields},
Discrete Analysis (2018), Paper No.~18, 42~pp.

\bibitem{phamV}
T.~Pham, L.~A.~Vinh, and F.~de Zeeuw,
\textit{Three-variable expanding polynomials and higher-dimensional
distinct distances},
Combinatorica \textbf{39}(2) (2019), 411--426.

\bibitem{quy}
M.-Q.~Pham,
\textit{On Falconer type functions and the distance set problem},
Mathematische Annalen \textbf{395}(1) (2026), Article~7.

\bibitem{sak}
A.~S\'ark\"ozy,
\textit{On sums and products of residues modulo \(p\)},
Acta Arithmetica \textbf{118}(4) (2005), 403--409.

\bibitem{shpas}
I.~E.~Shparlinski,
\textit{On the solvability of bilinear equations in finite fields},
Glasgow Mathematical Journal \textbf{50}(3) (2008), 523--529.

\bibitem{tao}
T.~Tao,
\textit{Expanding polynomials over finite fields of large characteristic,
and a regularity lemma for definable sets},
Contributions to Discrete Mathematics \textbf{10}(1) (2015), 22--98.

\bibitem{vinh}
L.~A.~Vinh,
\textit{On four-variable expanders in finite fields},
SIAM Journal on Discrete Mathematics \textbf{27}(4) (2013), 2038--2048.
\end{thebibliography}
\end{document}